\input amstex

\newif\ifproofmode                      
\proofmodefalse				

\newif\ifforwardreference		
\forwardreferencefalse			

\newif\ifchapternumbers			
\chapternumbersfalse			

\newif\ifcontinuousnumbering		
\continuousnumberingfalse		

\newif\iffigurechapternumbers		
\figurechapternumbersfalse		

\newif\ifcontinuousfigurenumbering	
\continuousfigurenumberingfalse		

\font\eqsixrm=cmr6			
\def\marginstyle{\eqsixrm}		

\newtoks\chapletter			
\newcount\chapno			
\newcount\eqlabelno			
\newcount\figureno			

\chapno=0
\eqlabelno=0
\figureno=0


\def\chapfolio{\ifnum \chapno>0 \the\chapno \else \the\chapletter \fi}


\def\bumpchapno{\ifnum \chapno>-1 \global \advance \chapno by 1
	\else \global \advance \chapno by -1 \setletter\chapno \fi
	\ifcontinuousnumbering \else \global\eqlabelno=0 \fi
	\ifcontinuousfigurenumbering \else \global\figureno=0 \fi}

%


%

\def\tempsetletter#1{\ifcase-#1 {}\or{} \or\chapletter={A}\or\chapletter={B}
  \or\chapletter={C} \or\chapletter={D} \or\chapletter={E}
  \or\chapletter={F} \or\chapletter={G} \or\chapletter={H}
  \or\chapletter={I} \or\chapletter={J} \or\chapletter={K}
  \or\chapletter={L} \or\chapletter={M} \or\chapletter={N}
  \or\chapletter={O} \or\chapletter={P} \or\chapletter={Q}
  \or\chapletter={R} \or\chapletter={S} \or\chapletter={T}
  \or\chapletter={U} \or\chapletter={V} \or\chapletter={W}
  \or\chapletter={X} \or\chapletter={Y} \or\chapletter={Z}\fi}

%

\def\chapshow#1{\ifnum #1>0 \relax #1%
   \else {\tempsetletter{\number#1}\chapno=#1 \chapfolio} \fi}

%
\def\today{\ifcase\month\or
January\or February\or March\or April\or May\or June\or
July\or August\or September\or October\or November\or December\fi
\space\number\day, \number\year}

%

\def\initialeqmacro{\ifproofmode
 \headline{\tenrm \today\hfill \jobname\ --- draft\hfill\folio}
     \hoffset=-1cm \immediate\openout2=allcrossreferfile \fi
 \ifforwardreference \input labelfile
     \ifproofmode \immediate\openout1=labelfile \fi \fi}


%

\def\chaplabel#1{\bumpchapno \ifproofmode \ifforwardreference
   \immediate\write1{\noexpand\expandafter\noexpand\def
   \noexpand\csname CHAPLABEL#1\endcsname{\the\chapno}}\fi\fi
   \global\expandafter\edef\csname CHAPLABEL#1\endcsname
   {\the\chapno}\ifproofmode\llap{\hbox{\marginstyle #1\ }}\fi\chapfolio}

%
\def\eqnum{\global\advance\eqlabelno by 1
   \eqno(\ifchapternumbers\chapfolio.\fi\the\eqlabelno)}

\def\eqlabel#1{\global\advance\eqlabelno by 1 \ifproofmode\ifforwardreference
 \immediate\write1{\noexpand\expandafter\noexpand\def
 \noexpand\csname EQLABEL#1\endcsname{\the\chapno.\the\eqlabelno?}}\fi\fi
 \global\expandafter\edef\csname EQLABEL#1\endcsname
 {\the\chapno.\the\eqlabelno?} \eqno(\ifchapternumbers\chapfolio.\fi
 \the\eqlabelno)\ifproofmode\rlap{\hbox{\marginstyle #1}}\fi}

\def\leqlabel#1{\global\advance\eqlabelno by 1 \ifproofmode\ifforwardreference
 \immediate\write1{\noexpand\expandafter\noexpand\def
 \noexpand\csname EQLABEL#1\endcsname{\the\chapno.\the\eqlabelno?}}\fi\fi
 \global\expandafter\edef\csname EQLABEL#1\endcsname
 {\the\chapno.\the\eqlabelno?} \leqno(\ifchapternumbers\chapfolio.\fi
 \the\eqlabelno)\ifproofmode\rlap{\hbox{\marginstyle #1}}\fi}

\def\eqalignnum{\global\advance\eqlabelno by 1
   &(\ifchapternumbers\chapfolio.\fi\the\eqlabelno)}

\def\eqalignlabel#1{\global\advance\eqlabelno by1 \ifproofmode
 \ifforwardreference\immediate\write1{\noexpand\expandafter\noexpand\def
 \noexpand\csname EQLABEL#1\endcsname
     {\the\chapno.\the\eqlabelno?}}\fi\fi
 \global\expandafter\edef\csname EQLABEL#1\endcsname
 {\the\chapno.\the\eqlabelno?}\ifchapternumbers\chapfolio.\fi
 \the\eqlabelno\ifproofmode\rlap{\hbox{\marginstyle
 #1}}\fi}

\def\eqref#1{(\ifundefined{EQLABEL#1}***\ifproofmode\ifforwardreference)%
   \else \write16{ ***Undefined Equation Reference #1*** }\fi
   \else \write16{ ***Undefined Equation Reference #1*** }\fi
   \else \edef\LABxx{\getlabel{EQLABEL#1}}%
   \def\LAByy{\expandafter\stripchap\LABxx}\ifchapternumbers
   \chapshow{\LAByy}.\expandafter\stripeq\LABxx
   \else\ifnum \number\LAByy=\chapno \relax\expandafter\stripeq\LABxx
   \else\chapshow{\LAByy}.\expandafter\stripeq\LABxx\fi\fi)\fi
   \ifproofmode\write2{Equation #1}\fi}

%

\def\fignum{\global\advance\figureno by 1 \relax
   \iffigurechapternumbers\chapfolio.\fi\the\figureno}\

\def\figlabel#1{\global\advance\figureno by 1\relax
 \ifproofmode\ifforwardreference
 \immediate\write1{\noexpand\expandafter\noexpand\def
 \noexpand\csname FIGLABEL#1\endcsname{\the\chapno.\the\figureno?}}\fi\fi
 \global\expandafter\edef\csname FIGLABEL#1\endcsname
 {\the\chapno.\the\figureno?}\iffigurechapternumbers\chapfolio.\fi
 \ifproofmode$^{\hbox{\marginstyle #1}}$\relax\fi\the\figureno}

\def\figref#1{\ifundefined{FIGLABEL#1}!!!!\ifproofmode\ifforwardreference)%
   \else \write16{ ***Undefined Equation Reference #1*** }\fi
   \else \write16{ ***Undefined Equation Reference #1*** }\fi
   \else \edef\LABxx{\getlabel{FIGLABEL#1}}%
   \def\LAByy{\expandafter\stripchap\LABxx}%
   \iffigurechapternumbers\chapshow{\LAByy}.\expandafter\stripeq\LABxx
   \else\ifnum\number\LAByy=\chapno \relax\expandafter\stripeq\LABxx
   \else\chapshow{\LAByy}.\expandafter\stripeq\LABxx\fi\fi
   \ifproofmode\write2{Figure #1}\fi\fi}

%

%

\def\getlabel#1{\csname#1\endcsname}
\def\ifundefined#1{\expandafter\ifx\csname#1\endcsname\relax}
\def\stripchap#1.#2?{#1}
\def\stripeq#1.#2?{#2}

\figurechapternumberstrue  

\chapternumberstrue        

\def\thmlbl#1{\figlabel{#1}}
\def\thmref{\figref}
\def\eqnlbl#1{\leqlabel{#1}}

\def\eqnref#1{\eqref{#1}}
\def\sectionnumber{\chapno}
\def\theoremnumber{\figureno}
\def\equationnumber{\eqlabelno}

\documentstyle{amsppt}
\overfullrule=0pt
\magnification =\magstep1
\baselineskip=18pt
\vcorrection{-.33truein}
\pageheight{9.0truein}
\input amssym.def

\def\BG{\bar G}
\def\tG{\tilde{\CalG}}

\def\BbbR{\Bbb{R}}

\def\CalG{\Cal{G}}

\def\CalL{\Cal{L}}
\def\CalR{\Cal{R}}

\def\CalU{\Cal{U}}
\def\CalV{\Cal{V}}

\def\R{\roman{Re\ }}

\def\sgn{\text{\rm sgn\ }}

\def\errfn{\text{\rm errfn\ }}

\def\deltdot{{\cdot \atop {\raise8pt\hbox{$\delta$}}}}
\def\myqed{\vrule height3pt depth2pt width3pt \bigskip}
\def\newsection{\centerline}


\pagewidth{6truein}
\pageheight{9truein}
\topmatter
\title
{Stability of viscous shock profiles for dissipative symmetric
hyperbolic-parabolic systems}
\endtitle
\leftheadtext{Stability of viscous shock profiles}
\rightheadtext{Corrado Mascia and Kevin Zumbrun}
\thanks
The second author thanks 
Instituto per le Applicazioni del Calcolo 
``M. Picone'' (CNR) and European TMR Project ``Hyperbolic Systems of 
Conservation Laws'' for their hospitality and for making possible
the visit (April 10--May 10, 2000) 
during which this work was carried out.
Research of the first author was supported in part by 
European TMR Project ``Hyperbolic Systems of Conservation Laws''.
Research of the second author was supported
in part by the National Science Foundation under Grants No. DMS-9107990
and DMS-0070765.
\endthanks
\abstract
Combining pointwise Green's function bounds obtained
in a companion paper [MZ.2] with earlier, spectral stability
results obtained in [HuZ],
we establish nonlinear orbital stability 
of small amplitude viscous shock profiles for the class of 
dissipative symmetric
hyperbolic-parabolic systems identified by Kawashima [Kaw],
notably including compressible Navier--Stokes
equations and the equations of magnetohydrodynamics,
obtaining sharp rates of decay in
$L^p$ with respect to small $L^1\cap H^3$ perturbations, 
$2\le p\le \infty$.
Our analysis follows the approach introduced in
[MZ.1] to treat stability of relaxation profiles.
\endabstract
\author
{Corrado Mascia and Kevin Zumbrun}
\endauthor
\date{August 3, 2001}
\enddate
\address
Dipartimento di Matematica ``G. Castelnuovo'',
Universit\`a di Roma ``La Sapienza'',
P.le Aldo Moro, 2 - 00185 Roma (ITALY)
\endaddress
\email
mascia\@ mat.uniroma1.it
\endemail 
\address
Department of Mathematics,
Indiana University,
Bloomington, IN  47405-4301
\endaddress
\email
kzumbrun\@indiana.edu
\endemail 
\endtopmatter
\document
\newsection {\bf Section 1. Introduction}
\sectionnumber=1 
\theoremnumber=0
\equationnumber=0
\smallskip
\TagsOnLeft

Consider the class of degenerate parabolic conservation laws
of {\it dissipative, symmetric hyperbolic--parabolic type} in 
the sense of Kawashima [Kaw], i.e., systems
$$
G(\tilde U)_t + F(\tilde U)_x = (B(\tilde U) \tilde U_x)_x,
\eqnlbl{general}
$$
satisfying
\medskip
(A1)\quad
$dF$, $dG$, $B$ symmetric, $dG>0$, $B\ge 0$ 
(symmetric hyperbolic--parabolicity),
\medskip
\medskip
(A2)\quad
No eigenvector of $dF dG^{-1}$ lies in the kernel of $BdG^{-1}$ (dissipativity),
\medskip
\noindent and
\medskip
(A3) \quad
$
\tilde U=\pmatrix \tilde u \\ \tilde v\endpmatrix$,
\quad
$
B=
\pmatrix 0 & 0 \\ 0 & b \endpmatrix
$,
\quad
$b>0$
(block structure),
\medskip
\noindent
in some neighborhood $\Cal{U}$ of a particular base point $\tilde U_*$,
where $\tilde U\in \BbbR^n$, 
$\tilde u\in  \BbbR^{n-r}$, $\tilde v\in \BbbR^r$, 
and $b\in \BbbR^{r\times r}$.
As discussed in [Kaw], this class of equations includes many physical
models arising in continuum mechanics: in particular, 
compressible Navier--Stokes equations and
the equations of compressible magnetohydrodynamics (MHD).
We make also the regularity assumption
\medskip
(H0) \quad $F$, $G$, $B\in C^3$,
\medskip
\noindent needed for our later analysis.

An interesting class of solutions of \eqnref{general} are 
{\it viscous shock profiles}, or asymptotically constant 
traveling wave solutions
$$
\tilde U=\bar U(x-st);
\quad 
\lim_{z\to \pm \infty} \bar U(z)=\tilde U_\pm
\eqnlbl{profile}
$$
connecting endstates $\tilde U_\pm$ corresponding
to (discontinuous) shock waves of the associated hyperbolic system
$$
G(\tilde U)_t + F(\tilde U)_x = 0.
\eqnlbl{hyperbolic}
$$
Existence of viscous profiles has been treated in the large
for compressible Navier--Stokes by, e.g., 
Gilbarg [Gi] 
and for MHD by Conley and Smoller [CS] (see [Fre.1] for a treatment
in the small amplitude case).
Most recently, Freist\"uhler [Fre.2] has carried out a general treatment
of existence in the small, 
analogous to that of Majda and Pego in the strictly parabolic case [MP],
for arbitrary models \eqnref{general} satisfying (A1)--(A3).

Let
$$
a_1(u)\le \dots \le a_n(u)
$$
denote the eigenvalues of $A_G:=dFdG^{-1}(\tilde U)$, 
$r_j(\tilde U)$ and $l_j(\tilde U)$ a smooth 
choice of associated right and left eigenvectors, $l_j\cdot r_k=\delta^j_k$,
and assume 
within neighborhood $\Cal{U}$ of base point $\tilde U_*$ that:

\medskip
(H1) \quad  The $p$th characteristic
field is of multiplicity one, i.e. $a_p(\tilde U)$ is a simple
eigenvalue of $A_G(\tilde U)$.

(H2) \quad  The $p$th characteristic field is genuinely
nonlinear, i.e. $\nabla a_p\cdot r_p (\tilde U)\ne 0$.
\medskip
(H3)\quad
The eigenvalues of $dF_{11}dG_{11}^{-1}$ 
(necessarily real and semisimple) are: (i) of constant
multiplicity; and, (ii) different from $a_p$.
\medskip
\noindent
Then, we have:

\proclaim{Proposition \thmlbl{freist1} ([Fre.2])}
Let (A1)--(A3), and (H1)--(H3) hold, and $F$, $G$, $B\in C^2$
in \eqnref{general} (implied in particular by (H0)).
Then, for left and right
states $\tilde U_\pm$ lying within a sufficently small neighborhood
$\CalV\subset \Cal{U}$ of $\tilde U_*$, and speeds $s$ lying within
a sufficiently small neighborhood of $a_p(\tilde U_*)$,
there exists a viscous profile \eqnref{profile} that is
``local'' in the sense that the image of $\bar u(\cdot)$ lies
entirely within $\Cal{V}$ if and only if the triple $(\tilde U_-,\tilde U_+,s)$
satisfies both the Rankine--Hugoniot relations: 
$$
s[G]=[F],
\tag{RH}
$$
and the Lax characteristic conditions for a $p$-shock:
$$
a_p(\tilde U_-)>s>a_p(\tilde U_+); \quad \sgn (a_j(\tilde U_-)-s)=\sgn (a_j(\tilde U_+)-s)\ne 0 \text{ for }
j\ne p.
\tag{L}
$$
(Note: The structure theorem of Lax [La,Sm] implies that (RH), always a 
{\it necessary} condition for existence of profiles, holds for 
$\tilde U_\pm\in \Cal{V}$ only if $s$ lies near some $a_j(\tilde U_*)$;
thus, the restriction on speed $s$ is only the assumption that the 
triple $(\tilde U_-,\tilde U_+,s)$ be associated with the $p$th and not some
other characteristic field).  
\endproclaim

{\bf Remark \thmlbl{ideal}.}
The generically satisfied conditions
(H1) and (H2) are implied by strict hyperbolicity
and genuine nonlinearity, respectively, of the associated 
hyperbolic system \eqnref{hyperbolic}.
In particular, they hold always for the equations of compressible 
gas dynamics, under the assumption of an ideal gas [Sm].
Condition (H3)(ii) corresponds to the requirement that ``hyperbolic,''
or unsmoothed modes in the solution be noncharacteristic, 
at least as regards the principal characteristic speed $a_p$.
This technical condition is needed in order that
the traveling wave ODE be of nondegenerate type, and is a standard
assumption in the theory.
See assumption (4) of [Fre.2], assumption ($\tilde{\text{\rm H}}$1)(ii) of [Z.3],
Appendix A.2, or assumption (H5) of [SZ] for restatements in
various different contexts.
Condition (H3)(i) is an additional technical assumption
that was used in the detailed Green's function analysis carried out
in [MZ.2] for slightly more general systems; we suspect that it can 
be dropped in the present, symmetrizable case.
Conditions (H3)(i)--(ii) are likewise satisfied for compressible 
ideal gas dynamics, for which $dF_{11}dG_{11}^{-1}$ is $1\times 1$
and identically equal to particle velocity ($u$ in Eulerian
coordinates, $0$ in Lagrangian coordinates).
\medskip

Stability of viscous profiles has been examined for compressible
Navier--Stokes equations in [MN,KMN,L.2], vith various partial results
concerning special (mainly zero-mass) initial data; 
the results of [KMN,L.2] are restricted to small amplitude profiles, while
the results of [MN] in the case of an isentropic $\gamma$-law gas
apply to profiles of amplitude $\alpha(\gamma)$, with 
$\alpha\to \infty$ as $\gamma \to 1$.
Yet, as pointed out recently in the fluid-dynamical survey [Te], the basic 
problem of stability of gas-dynamical shocks in its full generality 
remains open even for small amplitude waves, 
a significant gap in the theory of compressible gas dynamics.
More recently, Humpherys and Zumbrun [HuZ] have established {\it strong
spectral stability} of general, small amplitude shock waves of Kawashima
class systems, of the type constructed by Freist\"uhler, generalizing
a corresponding result of Goodman [Go.1--2] in the strictly parabolic 
case (see Section 2, below).
As discussed in [ZH,HuZ], strong spectral stability is roughly equivalent 
to, but slightly weaker than stability with respect to zero-mass perturbations.

In [MZ.2], applying the general machinery developed in [ZH,MZ.1], 
we have shown, for a slightly more general class of systems and for
profiles of arbitrary amplitude and type,
that strong spectral stability, plus hyperbolic stability of the corresponding
ideal shock (always satisfied for weak Lax shocks satisfying (H1)--(H3)),
are necessary and sufficient conditions for linearized orbital stability,
and, moreover, yield extremely detailed pointwise bounds on the Green's
function of the linearized evolution equations, analogous to those
obtained for relaxation profiles in [MZ.1].
In combination with the result of [HuZ], these results imply 
$L^1\cap L^p\to L^p$ {\it linearized orbital stability} of small-amplitude 
shock profiles of Kawashima class systems, with sharp rates of decay 
for all $1\le p\le \infty$.

The purpose of the present paper, extending the work of [HuZ,MZ.2],
is to establish $L^1\cap H^3\to L^p$ {\it nonlinear orbital stability} 
as solutions of \eqnref{general} of 
small amplitude profiles of general Kawashima class systems,
with sharp rates of decay for all $1\le p\le \infty$.
More precisely, we shall establish:

\proclaim{Theorem \thmlbl{nonlin}}
Let there hold (A1)--(A3) and (H0)--(H3) for
a general relaxation model \eqnref{general},
with $\tilde U_*$ and $\Cal{U}$ as above.
Then, for $\Cal{V}\subset \Cal{U}$ sufficiently small,
the viscous profiles $\Bar U$ described in Proposition 
\thmref{freist1} are nonlinearly orbitally stable
from $L^1\cap H^3$ to $L^p$, for all $p\ge 2$.
More precisely, for initial perturbations $U_0:=\tilde U_0-\bar U$
such that $|U_0|_{L^1\cap H^3}\le \zeta_0$, 
$\zeta_0$ sufficiently small,
the solution $\tilde U=(\tilde u,\tilde v) (x,t)$ of \eqnref{general} with
initial data $\tilde U_0$ satisfies
$$
|\tilde U(x,t)-\bar U(x-\delta(t))|_{L^p}\le C \zeta_0
(1+t)^{-\frac{1}{2}(1-1/p)}
\eqnlbl{3.31}
$$
and
$$
|\left(\tilde v(x,t)-\bar v(x-\delta(t))\right)_x|_{L^p}\le C \zeta_0
(1+t)^{-\frac{1}{2}(1-1/p)},
\eqnlbl{deriv3.31}
$$
for all $2\le p\le \infty$, for some $\delta(t)$ satisfying
$$
|\dot \delta (t)|\le C \zeta_0 (1+t)^{-\frac{1}{2}}
\eqnlbl{3.32}
$$
and
$$
|\delta(t)|\le C \zeta_0.
\eqnlbl{3.33}
$$
\endproclaim

\noindent In particular, 
this implies nonlinear stability
of all small amplitude gas-dynamical profiles under the
assumptions of an ideal gas; 
see Remark \thmref{ideal} above. 

Our analysis in this paper, and in the companion paper [MZ.2],
follows closely to that used in [MZ.1] to
treat stability of relaxation profiles, making use of
structural similarities between the two problems 
pointed out in [Ze.2,Z.3];
in turn, the argument of [MZ.1] makes substantial use of the
pointwise semigroup machinery introduced in [ZH,Z.2] to 
treat the parabolic case.
As are general features of the pointwise semigroup approach, 
we obtain through this program extremely detailed
pointwise bounds on the Green's function (more properly,
distribution) of the linearized operator about the wave,
sufficiently strong that the nonlinear stability analysis
becomes comparatively simple.
Moreover, these bounds (though not the nonlinear stability argument
in which they are used, see Remark 1.5 just below) 
depend only on a generalized spectral stability 
(i.e., Evans function) condition and not on the amplitude or 
type (Lax, under-, or overcompressive) of the wave.

{\bf Remark \thmlbl{ngnl}.} The assumption of genuine nonlinearity (H2)
is not needed either for the existence or the stability result,
but is made only to simplify the discussion.
Though we stated above only the restriction
to the genuinely nonlinear case,
existence was in fact treated for the general (nongenuinely nonlinear) 
case in [Fre.2], substituting for the Lax entropy condition (L)
the strict Liu entropy inequality (E) of [L.3].
Likewise, as described in [HuZ], the result of spectral stability
may be extended to the general case by substituting for the 
``Goodman-type'' weighted energy estimates of [HuZ], Section 5, 
the variation introduced by Fries [Fri.1--2]
to treat the nongenuinely nonlinear case for strictly parabolic viscosities.
Accordingly, the result of Theorem \thmref{nonlin} holds also in this
case.
\medskip
{\bf Remark \thmlbl{gammalaw}.} The restriction to small-amplitude
shocks arises only through the energy estimates used to close the
nonlinear iteration argument.
In particular, it should be possible for the isentropic, $\gamma$-law
gas case to treat amplitudes of the same order $\alpha(\gamma)$ treated in [MN].
A fundamental open problem is to remove the amplitude
restriction altogether, as done in the
strictly parabolic case [ZH,Z.2] and in the case of discrete kinetic
relaxation models [MZ.1], replacing it with 
the generalized spectral condition of the linearized theory.
\medskip

{\bf Remark \thmlbl{nec}.}
Useful necessary conditions for viscous stability have been obtained
in [Z.3] for arbitrary amplitude profiles of the more general class 
of models considered in [MZ.2], using the stability index of [GZ,BSZ].
Strengthened versions of the one-dimensional inviscid stability criteria
of Erpenbeck--Majda [Er,M.1--3], these
readily yield examples of unstable large-amplitude profiles,
similarly as in the strictly parabolic case (see, 
e.g., [GZ,FreZ,ZS,Z.3]).
This shows that the spectral stability requirement is not vacuous
in the large amplitude case.

\medskip
{\bf Plan of the paper.}
In Section 2 we cite the spectral stability result of [HuZ],
and in Section 3 the pointwise Green's function bounds of [MZ.2].
As an immediate corollary, we establish in Section 4 the 
linearized orbital stability of general shock profiles satisfying 
the necessary conditions of spectral and hyperbolic stability:
in particular, of small amplitude profiles of Kawashima class systems.
Finally, we establish in Section 5 the main result of
nonlinear stability of weak profiles of Kawashima class systems,
by a modified version of the argument of [MZ.1].

\bigskip
{\bf Note}:  
Liu and Zeng have informed us [LZe.2] that they also have obtained 
nonlinear stability of weak Navier--Stokes profiles,
by a different argument based on the approximate 
Green's function approach of [L.2].
This method is inherently limited to weak shocks of classical, Lax type,
both at the linear and nonlinear level.
\bigskip

\newsection {\bf Section 2. Spectral stability.}
\sectionnumber=2 
\theoremnumber=0
\equationnumber=0
\smallskip
\TagsOnLeft
\bigskip

Take without loss of generality $s=0$, so that
$\tilde U=\bar U(x)$ becomes a stationary solution.
Then, the linearized equations of \eqnref{general} about 
$\bar U$ take the form
$$
U_t =LU:= -(A^0)^{-1}(AU)_x + (A^0)^{-1}(BU_x)_x,
\eqnlbl{linearized}
$$
where 
$$
B:= B(\bar U), 
\quad A^0:=dG(\bar U), 
\quad  Av:= dF(\bar U)v
- (dB v)\bar U_x.
,
\eqnlbl{AB}
$$

{\bf Definition \thmlbl{stable}.}  We call the profile $\bar U(\cdot)$
{\it strongly spectrally stable} if the linearized operator
$L$ about the wave has no spectrum in the closed unstable
complex half-plane $\{\lambda: \, \text{\rm Re }\lambda \ge 0\}$
except at the origin, $\lambda=0$.
(Recall, [Sat], that $\lambda=0$ is always in the spectrum of $L$,
since $L\bar U_x=0$ by direct calculation/differentiation of
the traveling wave ODE).
\medskip

The spectral stability of small-amplitude viscous profiles of
the slightly more more general class of dissipative, {\it symmetrizable}
degenerate parabolic systems has been investigated in [HuZ].
To apply these results, we have only to note that \eqnref{general}, 
expressed with respect to variable $\tilde G:=G(\tilde U)$, becomes
$$
\tilde G_t + F(\tilde U(\tilde G))_x= (B(\tilde U(\tilde G))\tilde U(\tilde G)_x)_x,
\eqnlbl{symmetrizable1}
$$
or, in quasilinear form,
$$
\tilde G_t + F_{\tilde U}G_{\tilde U}^{-1}\tilde G_x = (B G_{\tilde U}^{-1} \tilde G_x)_x.
\eqnlbl{symmetrizable}
$$
System \eqnref{symmetrizable} is evidently symmetrizable
(by the symmetric positive definite $G_{\tilde U}^{-1}$), and
{\it dissipative} (since the dissipativity condition is independent
of coordinate system), with the additional {\it block structure}
property that the left kernel of the new viscosity matrix
$B G_{\tilde U}^{-1}$ is constant.
This is almost the class of equations considered in [HuZ], the
difference being that there the block structure assumption was
that the viscosity matrix have constant right instead of left
kernel.  This discrepancy is unimportant in the analysis, since,
at the linearized level, one  case may be converted to the
other by the (linear) change of coordinates
$A^0 V:=U$, corresponding to similarity
transform $M\to (A^{0})^{-1} M A^0$.
Alternatively, in the particular case of our interest, we may simply
carry out energy estimates in the natural coordinates of \eqnref{linearized};
for related calculations, see Section 4.1, below.

Thus, the results of [HuZ] apply to the somewhat larger class
of dissipative symmetrizable systems of form 
\eqnref{symmetrizable1}--\eqnref{symmetrizable1},
with the block structure condition that {\it either} the left or
right kernel of the viscosity matrix be constant, and we may
conclude, in particular:

\proclaim{Theorem \thmlbl{HuZ}[HuZ]}
Let (A1)--(A3) and (H0)--(H3) hold,
and let $\bar U(x-st)$ be a viscous shock solution of \eqnref{general}
such that the profile $\{\bar u(z)\}$ lies entirely within a sufficently 
small neighborhood $\CalV\subset \Cal{U}$ of $\tilde U_*$, and the speed 
$s$ lies within a sufficiently small neighborhood of $a_p(\tilde U_*)$.
Then, $\bar U$ is strongly spectrally stable, in the sense of Definition
\thmref{stable} above.
\endproclaim
\bigskip

\newsection {\bf Section 3. Pointwise Green's function bounds.}
\sectionnumber=3
\theoremnumber=0
\equationnumber=0
\smallskip
\TagsOnLeft
\bigskip

In a companion paper to this one [MZ.1], we have investigated linearized
stability and behavior of viscous profiles of systems of the general
form \eqnref{symmetrizable}, not necessarily symmetrizable, satisfying
the standard set of conditions identified in [Z.3].
As remarked in [Z.3,MZ.1], these hold always for symmetric systems \eqnref{general} 
satisfying (A1)--(A3), (H0)--(H3), and the small-amplitude profiles described in 
Proposition \thmref{freist1};
however, they may also hold in much greater generality, in particular for 
shocks of large amplitude, or nonclassical type.
For simplicity of exposition, we shall restrict our discussion here to the 
present case of interest; for the general case, see [MZ.1].

The first main result of [MZ.2], generalizing the corresponding
results established for viscous, strictly parabolic, shocks in [ZH],
and for relaxation shocks in [MZ.1], is:

\proclaim{Theorem \thmlbl{D}[MZ.2]}
Under assumptions (A1)--(A3), (H0)--(H3), 
small amplitude shock profiles are $L^1\cap L^p\to L^p$
linearly orbitally stable for $p>1$ if and only if
they are strongly spectrally stable, with sharp decay bounds
$$
|U(\cdot, t)+\delta(t)\bar U'(\cdot)|_{L^p}\le 
C(1+t)^{-\frac{1}{2}(1-1/p)}(|U_0|_{L^1}+ |U_0|_{L^p})
\eqnlbl{lindecaybound1}
$$
for initial data $U_0\in L^1\cap L^p$, some choice of $\delta$.
\endproclaim

Theorem \thmref{D} is obtained (see, e.g., Corollary 4.2, below for
one direction) as a consequence of
detailed, pointwise bounds on the Green's function 
(more properly speaking, distribution)
$\CalG(x,t;y)$ of the linearized evolution equations 
$$
G_t= L_G G
:=- (A_GG)_x + (B_G G_x)_x,
\eqnlbl{linsymmetrizable}
$$
$$
B_G:= B (A^0)^{-1}(\bar U(x)),
\quad
A_Gv:= dF (A^0)^{-1} v - dB(\bar U(x))v \bar U'(x),
\eqnlbl{Gcoeffs}
$$
corresponding to \eqnref{symmetrizable}, i.e., expressed
with respect to the conservative variable $G$.
We now describe these bounds, for use in the following sections.

Let $a_j^\pm$, $j=1,\dots (n)$ denote the eigenvalues of
$A_G(\pm\infty)$, and $l_j^\pm$ and $r_j^\pm$ 
associated left and right eigenvectors, respectively,
normalized so that $l_j^{\pm t} r_k^{\pm}=\delta^j_k$. 
In case $A_G^\pm$ is strictly hyperbolic, these are uniquely defined.
In the general case, we require further that $l_j^\pm$, $r_j^\pm$
be left and right eigenvectors also of $P_j^{\pm } B_G^\pm P_j^\pm$, 
$P_j^\pm:=R_j^\pm L_j^{\pm t}$, where
$L_j^\pm$ and $R_j^\pm$ denote $m_j^\pm\times m_j^\pm$ 
left and right eigenblocks associated
with the $m_j^\pm$-fold eigenvalue $a_j^\pm$, normalized so that 
$L_j^\pm R_j^\pm= I_{m_j^\pm}$.
(Note: The matrix $P_j^{\pm } B_G^\pm P_j^\pm\sim
L_j^{\pm t} B_G^\pm R_j^\pm$ 
is necessarily diagonalizable,
by simultaneous symmetrizability of $A_G$, $B_G$).

Eigenvalues $a_j(x)$, and eigenvectors $l_j$, $r_j$
correspond to large-time convection rates and modes of propagation 
of the degenerate model \eqnref{symmetrizable}.
Likewise, let $a^{*}_j(x)$, $j=1,\dots,(n-r)$ 
denote the eigenvalues of 
$$
\aligned
A^*_G&:= A_{G,11}- A_{G,12} B_{G,22}^{-1}B_{G,21}\\
&=
A_{11}(A^0_{11})^{-1} ,
\endaligned
\eqnlbl{A*}
$$
and $l^*_j(x)$, $r^*_j(x)\in \BbbR^{n-r}$ 
associated left and right eigenvectors,
normalized so that $l^{*t}_jr_j\equiv \delta^j_k$.
More generally, for an $m_j^*$-fold eigenvalue, we choose
$(n-r)\times m_j^* $ blocks $L_j^*$ and $R_j^*$ of eigenvectors
satisfying the {\it dynamical normalization}
$$
L_j^{*t}\partial_x R_j^{*}\equiv 0,
$$
along with the usual static normalization 
$L^{*t}_jR_j\equiv \delta^j_kI_{m_j^*}$; as shown in Lemma 4.9, [MZ.1],
this may always be achieved with bounded $L_j^*$, $R_j^*$.
Associated with $L_j^*$, $R_j^*$, define extended, $n\times m_j^*$ blocks
$$
\CalL_j^*:=\pmatrix L_j^* \\ 0\endpmatrix,
\quad
\CalR_j^*:=
\pmatrix R_j^*\\
-B_{G,22}^{-1}B_{G,21} R_j^*\endpmatrix.
\eqnlbl{CalLR}
$$
%
Eigenvalues $a_j^*$ and eigenmodes $\CalL_j^*$, $\CalR_j^*$
correspond, respectively, to short-time hyperbolic characteristic speeds 
and modes of propagation for the reduced, hyperbolic part of degenerate
system \eqnref{symmetrizable}.
Note that our discussion in the introduction of condition (H3) is
validated by the second equality in \eqnref{A*}, a 
nontrivial consequence of symmetry/block structure in the 
original system \eqnref{general}.

Define time-asymptotic, {\it scalar diffusion rates}
$$
\beta_j^{\pm}:= \left(l_j^{t} B_G r_j\right)_\pm, \quad j=1,\dots, n,
\eqnlbl{beta}
$$
and local, $m_j\times m_j$ {\it dissipation coefficients}
$$
\eta_j^*(x):= -l_j^{*t} D^* r_j^* (x), \quad j=1,\dots,J\le n-r, 
\eqnlbl{eta}
$$
where
$$
{D^*}:= A_{G,12} \Big[A_{G,21}-A_{G,22}
B_{G,22}^{-1} B_{G,21}+ A^*_G 
B_{G,22}^{-1} B_{G,21} \Big],
\eqnlbl{D*}
$$
is an effective dissipation precisely analogous to the
effective diffusion predicted by formal, Chapman--Enskog expansion
in the (dual) relaxation case.
As described in Appendix A2 of [MZ.2], these quantities arise in a natural 
way, through Taylor expansion of the (frozen-coefficient) Fourier symbol 
$$
-i\xi A_G(x)-\xi^2 B_G(x)
\eqnlbl{frozenx}
$$
of the linearized operator $L_G$ about $\xi=0$ and $\xi=\infty$, respectively.

The important quantities $\eta_j^*$, $\beta_j$ were identified by
Zeng [Ze.1,LZe.1] in her study by Fourier transform techniques
of decay to {\it constant} 
(necessarily equilibrium) {\it solutions} 
$(\bar u, \bar v) \equiv (u_\pm,v_\pm)$ 
of relaxation systems, corresponding at the linearized level
to the study of the limiting equations 
$$
U_t=L_{G}^\pm U:= -A_G^{\pm} U_x+ B_G^\pm U_{xx}
\eqnlbl{limiting}
$$
as $x\to\pm \infty$ of the linearized evolution equations 
\eqnref{linsymmetrizable}.
As a consequence of dissipativity, (A2), we have (see, e.g., [Kaw,LZe.1]) 
that
$$
\beta_j^{\pm}>0, \quad 
j=1,\dots,n
\eqnlbl{goodbeta}
$$
and
$$
\R \sigma(\eta_j^*(x))>0, \quad
j=1,\dots,J\le n-r.
\eqnlbl{goodeta}
$$

\proclaim{Proposition \thmlbl{greenbounds}[MZ.2]}
For weak shock profiles, under assumptions (A1)--(A3), (H0)--(H3),
the Green's function $\CalG(x,t;y)$ associated with
the linearized evolution equations \eqnref{linsymmetrizable} 
may in the Lax or overcompressive case be decomposed as
$$
\CalG(x,t;y)= H + E+  S + R,
\eqnlbl{ourdecomp}
$$
where, for $y\le 0$:
$$
\aligned
H(x,t;y)&:=
\sum_{j=1}^{J} 
\CalR_j^*(x) \zeta_j^*(y,t)
\delta_{x-\bar a_j^* t}(-y)
\CalL_j^{*t}(y)\\
&=
\sum_{j=1}^{J} \CalR_j^*(x) \Cal{O}(e^{-\eta_0 t})
\delta_{x-\bar a_j^* t}(-y)
\CalL_j^{*t}(y),
\endaligned
\eqnlbl{multH}
$$
$$
\aligned
E&(x,t;y):=
\sum_{a_k^- > 0}
[c^{0}_{k,-}]
\bar G'(x)
l_k^{-t}
\left(\errfn\left(\frac{y+a_k^{-}t}{\sqrt{4\beta_k^{-}t}}\right)
-\errfn \left(\frac{y-a_k^{-}t}{\sqrt{4\beta_k^{-}t}}\right)\right),
\endaligned
\eqnlbl{E}
$$
and
$$
\aligned
S(x,t;y)&:=
\chi_{\{t\ge 1\}} 
\sum_{a_k^{-}<0}r_k^{-}  {l_k^{-}}^t
(4\pi \beta_k^-t)^{-1/2} e^{-(x-y-a_k^{-}t)^2 / 4\beta_k^{-}t} 
\\
&+ 
\chi_{\{t\ge 1\}} 
\sum_{a_k^{-} > 0} r_k^{-}  {l_k^{-}}^t
(4\pi \beta_k^{-}t)^{-1/2} e^{-(x-y-a_k^{-}t)^2 / 4\beta_k^{-}t}
\left({e^x \over e^x+e^{-x}}\right)\\
&+ 
\chi_{\{t\ge 1\}}
\sum_{a_k^{-} > 0, \,  a_j^{-} < 0} 
[c^{j,-}_{k,-}]r_j^{-}  {l_k^{-}}^t
(4\pi \bar\beta_{jk}^{-} t)^{-1/2} e^{-(x-z_{jk}^{-})^2 / 
4\bar\beta_{jk}^{-} t} 
\left({e^{ -x} \over e^x+e^{-x}}\right),\\
&+ 
\chi_{\{t\ge 1\}}
\sum_{a_k^{-} > 0, \,  a_j^{+} > 0} 
[c^{j,+}_{k,-}]r_j^{+}  {l_k^{-}}^t
(4\pi \bar\beta_{jk}^{+} t)^{-1/2} e^{-(x-z_{jk}^{+})^2 / 
4\bar\beta_{jk}^{+} t} 
\left({e^{ x} \over e^x+e^{-x}}\right)\\
\endaligned
\eqnlbl{S}
$$
denote hyperbolic, excited, and scattering terms, respectively, 
and $R$ denotes a faster decaying residual, satisfying: 
$$
\aligned
R(x,t;y)&= 
\Cal{O}(e^{-\eta(|x-y|+t)})\\
&+\sum_{k=1}^n 
\Cal{O} \left( (t+1)^{-1/2} e^{-\eta x^+} 
+e^{-\eta|x|} \right) 
t^{-1/2}e^{-(x-y-a_k^{-} t)^2/Mt} \\
&+
\sum_{a_k^{-} > 0, \, a_j^{-} < 0} 
\chi_{\{ |a_k^{-} t|\ge |y| \}}
\Cal{O} ((t+1)^{-1/2} t^{-1/2})
e^{-(x-a_j^{-}(t-|y/a_k^{-}|))^2/Mt}
e^{-\eta x^+}, \\
&+
\sum_{a_k^{-} > 0, \, a_j^{+}> 0} 
\chi_{\{ |a_k^{-} t|\ge |y| \}}
\Cal{O} ((t+1)^{-1/2} t^{-1/2})
e^{-(x-a_j^{+} (t-|y/a_k^{-}|))^2/Mt}
e^{-\eta x^-}, \\
\endaligned
\eqnlbl{Rbounds}
$$
$$
\aligned
R_y(x,t;y)&= 
\sum_{j=1}^J \Cal{O}(e^{-\eta t})\delta_{x-\bar a_j^* t}(-y) 
+
\Cal{O}(e^{-\eta(|x-y|+t)})\\
&+\sum_{k=1}^n 
\Cal{O} \left( (t+1)^{-1/2} e^{-\eta x^+} 
+e^{-\eta|x|} \right) 
t^{-1}
e^{-(x-y-a_k^{-} t)^2/Mt} \\
&+
\sum_{a_k^{-} > 0, \, a_j^{-} < 0} 
\chi_{\{ |a_k^{-} t|\ge |y| \}}
\Cal{O} ((t+1)^{-1/2} t^{-1}) 
e^{-(x-a_j^{-}(t-|y/a_k^{-}|))^2/Mt}
e^{-\eta x^+} \\
&+
\sum_{a_k^{-} > 0, \, a_j^{+} > 0} 
\chi_{\{ |a_k^{-} t|\ge |y| \}}
\Cal{O} ((t+1)^{-1/2} t^{-1}) 
e^{-(x-a_j^{+}(t-|y/a_k^{-}|))^2/Mt}
e^{-\eta x^-}, \\
\endaligned
\eqnlbl{Rybounds}
$$
$$
\aligned
R_x(x,t;y)&= 
\sum_{j=1}^J \Cal{O}(e^{-\eta t})\delta_{x-\bar a_j^* t}(-y) 
+
\Cal{O}(e^{-\eta(|x-y|+t)})\\
&+\sum_{k=1}^n 
\Cal {O} \left( (t+1)^{-1} e^{-\eta x^+} 
+e^{-\eta|x|} \right) 
t^{-1} (t+1)^{1/2}
e^{-(x-y-a_k^{-} t)^2/Mt} \\
&+
\sum_{a_k^{-} > 0, \, a_j^{-} < 0} 
\chi_{\{ |a_k^{-} t|\ge |y| \}}
\Cal{O}(t+1)^{-1/2} t^{-1}) 
e^{-(x-a_j^{-}(t-|y/a_k^-|))^2/Mt}
e^{-\eta x^+} \\
&+
\sum_{a_k^{-} > 0, \, a_j^{+} > 0} 
\chi_{\{ |a_k^{-} t|\ge |y| \}}
\Cal{O}(t+1)^{-1/2} t^{-1}) 
e^{-(x-a_j^{+}(t-|y/a_k^{-}|))^2/Mt}
e^{-\eta x^-}, \\
\endaligned
\eqnlbl{Rxbounds}
$$
and
$$
\aligned
R_{xy}(x,t;y)&= 
(\partial/\partial y)\big(\sum_{j=1}^J 
\Cal{O}(e^{-\eta t})\delta_{x-\bar a_j^* t}(-y) \big)\\
& +
\sum_{j=1}^J \Cal{O}(e^{-\eta t})\delta_{x-\bar a_j^* t}(-y) 
+
\Cal{O}(e^{-\eta(|x-y|+t)})\\
&+\sum_{k=1}^n 
\Cal {O} \left( (t+1)^{-3/2} e^{-\eta x^+} 
+e^{-\eta|x|} \right) 
t^{-3/2} (t+1)
e^{-(x-y-a_k^{-} t)^2/Mt} \\
&+
\sum_{a_k^{-} > 0, \, a_j^{-} < 0} 
\chi_{\{ |a_k^{-} t|\ge |y| \}}
\Cal{O}(t+1)^{-1/2} t^{-3/2}) 
e^{-(x-a_j^{-}(t-|y/a_k^-|))^2/Mt}
e^{-\eta x^+} \\
&+
\sum_{a_k^{-} > 0, \, a_j^{+} > 0} 
\chi_{\{ |a_k^{-} t|\ge |y| \}}
\Cal{O}(t+1)^{-1/2} t^{-3/2}) 
e^{-(x-a_j^{+}(t-|y/a_k^{-}|))^2/Mt}
e^{-\eta x^-}, \\
\endaligned
\eqnlbl{Rxybounds}
$$
for some $\eta$, $M>0$, where $x^\pm$ denotes the positive/negative 
part of $x$,  and indicator function $\chi_{\{ |a_k^{-}t|\ge |y| \}}$ is 
one for $|a_k^{-}t|\ge |y|$ and zero otherwise.
Symmetric bounds hold for $y\ge0$.

Here, the averaged convection rates
$ \bar a_j^*(y,t)$ in \eqnref{multH} and \eqnref{Rxbounds}--\eqnref{Rxybounds}
denote the time-averages over $[0,t]$ 
of $a_j^*(x)$ along characteristic paths $z_j^*=z_j^*(y,t)$ defined by
$$
dz_j^*/dt= a_j^*(z_j^*), \quad z_j^*(0)=y,
\eqnlbl{char}
$$
and the dissipation matrix 
$\zeta_j^*=\zeta_j^*(y,t)\in \BbbR^{m_j\times m_j}$ 
is defined by the {\it dissipative flow} 
$$
d\zeta_j^*/dt= -\eta_j^*(z_j^*)\zeta_j^*, \quad \zeta_j^*(y)=I_{m_j}.
\eqnlbl{diss}
$$
Similarly, in \eqnref{S},
$$
z_{jk}^{\pm(y,t)}:=a_j^{\pm}\left(t-\frac{|y|}{|a_k^{-}|}\right)
\eqnlbl{zjk}
$$
and
$$
\bar \beta^{\pm}_{jk}(x,t;y):= 
\frac{|x^\pm|}{|a_j^{\pm} t|} \beta_j^{\pm}
+
\frac{|y|}{|a_k^{-} t|} 
\left( \frac{a_j^\pm}{a_k^{-}}\right)^2 \beta_k^{-},
\eqnlbl{barbeta}
$$
represent, respectively, approximate scattered characteristic
paths and the time-averaged diffusion rates along those paths.
In all equations, $a_j$, $a_j^{*\pm}$, $l_j$, $\CalL_j^{*\pm }$,
$r_j$, $\CalR_j^{*\pm}$, 
$\beta_j^{\pm}$ and $\eta_j^*$ are as defined just above, 
and scattering coefficients $[c_{k,-}^{j,i}]$, $i=-,0,+$, are constants, 
uniquely determined by
$$
\sum_{a_j^{-} < 0} [c_{k, \, -}^{j, \, -}]r_j^{-} +
\sum_{a_j^{+} > 0} [c_{k, \, -}^{j, \, +}]r_j^{+} +
 [c_{k,-}^{0}] \big (G(+\infty)-G(-\infty) \big)
= r_k^{-}
\eqnlbl{scattering}
$$
for each $k=1,\dots n$, and satisfying
$$
\sum_{a_k^->0} [c_{k,-}^{0}] l_k^{-}
=
 \sum_{a_k^+<0} [c_{k,+}^{0}] l_k^{+} 
=\pi,
\eqnlbl{pi}
$$
where the constant vector $\pi$
is the left zero effective eigenfunction of $L_G$ associated
with the right eigenfunction $\bar G'$.

\endproclaim

Proposition \thmref{greenbounds}, the variable-coefficient
generalization of the constant-coefficient results of
[Ze.1,LZe.1],  was established in [MZ.2] by Laplace transform 
(i.e., semigroup) techniques generalizing the Fourier transform 
approach of [Ze.1--2,LZe];
for discussion/geometric interpretation, see [Z.2,MZ.1--2].
In our stability analysis, we will use only a small part of 
the detailed information given in the proposition, namely
$L^p\to L^q$ estimates on the time-decaying portion 
$H+S+R$ of the Green's function $\CalG$ (see Lemma 4.1, below).
However, the stationary portion $E$ of the Green's function
must be estimated accurately for an efficient stability analysis.

\proclaim{Lemma \thmlbl{goodproj}}
Under the assumptions of Proposition \thmref{greenbounds}, there
holds, additionally,
$$
\Pi_2 (A^0(x))^{-1}H(x,t;y)\equiv 0,
\eqnlbl{good1}
$$
where $\Pi_2:=\pmatrix 0&0\\
0&I_r\endpmatrix$
denotes projection onto the final $r$ coordinates of $G$.
\endproclaim

{\bf Proof.}
Equivalently, we must show that $\Pi_2 (A^0)^{-1}\CalR_j^*\equiv 0$
for all $1\le j\le J$.
This is most easily verified by the intrinsic property 
of $\CalL_j^*$ and $\CalR_j^*$ (readily seen from our formulae) 
that they lie, 
respectively, in the left and right kernel of $B_G$.
For, this gives
$$
\aligned
0\equiv B_G \CalR_j^*
&:= \pmatrix 0 & 0\\0& b\endpmatrix 
(A^0)^{-1}\CalR_j\\
&=
b \Pi_2 (A^0)^{-1}\CalR_j^*,
\endaligned
\eqnlbl{calc}
$$
yielding the result by invertibility of $b$, (see condition (A3)).
\myqed

Lemma \thmref{goodproj} quantifies the observation that, in $U=(u,v)^t$ 
coordinates, the ``parabolic'' variable $v$ experiences smoothing
under the evolution of \eqnref{linearized}, whereas the
``hyperbolic'' coordinate $v$ does not.
(Recall that $U=(A^0)^{-1}G$).

\bigskip

\newsection {\bf Section 4. Linearized stability.}
\sectionnumber=4 
\theoremnumber=0
\equationnumber=0
\smallskip
\TagsOnLeft
\bigskip

We now show, for the Lax case under consideration,
that linearized orbital stability 
follows immediately from the pointwise bounds 
of Proposition \thmref{greenbounds}, thus partially recovering
the result of [MZ.2] that was stated in Proposition \thmref{D}.
This analysis motivates the nonlinear argument to follow
in Section 5.

We carry out our analysis with respect to the
conservative variable $G=A^0(x)U$, i.e. with respect to 
linearized equations \eqnref{linsymmetrizable}.
Similarly as in [Z.2,MZ.1--2], 
define the {\it linear instantaneous projection}:
$$
\aligned
\varphi(x,t)&:=
\int_{-\infty}^{+\infty} E(x,t;y)G_0(y)\, dy
\\
&=:-\delta(t) \bar G'(x),\\
\endaligned
\eqnlbl{linproj}
$$
where $G_0$ denotes the initial data for \eqnref{linearized},
and $\bar G=G(\bar U(x))$ as usual.
The amplitude $\delta$ may be expressed, alternatively, as
$$
\delta(t)= 
-\int_{-\infty}^{+\infty} e(y,t)G_0(y)\, dy,
$$
where
$$
E(x,t;y)=:\bar G'(x)
e(y,t),
\eqnlbl{eE}
$$
i.e., 
$$
e(y,t):=
\sum_{a_k^{-}}
\left(\errfn\left(\frac{y+a_k^{-}t}{\sqrt{4\beta_k^{-}t}}\right)
-\errfn \left(\frac{y-a_k^{-}t}{\sqrt{4\beta_k^{-}t}}\right)\right)
l_k^{-}
\eqnlbl{e}
$$
for $y\le 0$, and symmetrically for $y\ge 0$.

Then, the solution $G$ of \eqnref{linearized} satisfies
$$
G(x,t)-\varphi(x,t)=
\int_{-\infty}^{+\infty}(H + \tG)(x,t;0)G_0(y)\, dy,
\eqnlbl{convolution}
$$
where 
$$
\tG:=S+R
\eqnlbl{tG}
$$
is the regular part and $H$ the singular part of
the time-decaying portion of the Green's function $\CalG$.

\proclaim{Lemma \thmlbl{2.05}}  
For weak shock profiles of \eqnref{general}, 
under assumptions (A1)--(A3), (H0)--(H3), there hold:
$$
|\int_{-\infty}^{+\infty} \tG(\cdot,t;y)f(y)dy|_{L^p}
\le C (1+t)^{-\frac{1}{2}(1-1/r)} |f|_{L^q},
\eqnlbl{tGbounds}
$$
$$
|\int_{-\infty}^{+\infty} \tG_y(\cdot,t;y)f(y)dy|_{L^p}
\le C (1+t)^{-\frac{1}{2}(1-1/r)-1/2} |f|_{L^q}
+ Ce^{-\eta t} |f|_{L^p},
\eqnlbl{tGybounds}
$$
$$
|\int_{-\infty}^{+\infty} 
\tG_x(\cdot,t;y)
f(y)dy|_{L^p}
\le C (1+t)^{-\frac{1}{2}(1-1/r)} |f|_{L^q}
+ Ce^{-\eta t} |f|_{L^p},
\eqnlbl{tGxbounds}
$$
$$
|\int_{-\infty}^{+\infty} 
\tG_{xy}(\cdot,t;y)
f(y)dy|_{L^p}
\le C (1+t)^{-\frac{1}{2}(1-1/r)-1/2} |f|_{L^q}
+ Ce^{-\eta t} |f|_{W^{1,p}},
\eqnlbl{tGxybounds}
$$
and
$$
|\int_{-\infty}^{+\infty} H(\cdot,t;y)f(y)dy|_{L^p}
\le Ce^{-\eta t} |f|_{L^p},
\eqnlbl{Hbounds}
$$
for all $t\ge 0$, some $C$, $\eta>0$, for any
$1\le q\le p$ 
(equivalently, $1\le r\le p$)
and $f\in L^q \cap W^{1,p}$, where $1/r+1/q=1+1/p$.
\endproclaim

{\bf Proof.}  
Bounds \eqnref{tGbounds}--\eqnref{tGxybounds} follow by
the Hausdorff-Young inequality together
with bounds \eqnref{S} and \eqnref{Rbounds}--\eqnref{Rxybounds},
precisely as in [Z.2,MZ.1--2].
Bound \eqnref{Hbounds} follows by direct computation and
the fact that particle paths $ z_j(y,t)$
satisfy uniform bounds 
$$
1/C \le |(\partial/\partial y)z_j| < C,
$$
for all $y$, $t$, by the fact that characteristic speeds $a_j(x)$ converge
exponentially as $x\to \pm \infty$ to constant states.
\myqed

\proclaim{Corollary \thmlbl{suff}}
Let $\bar U$ be a weak shock profile of \eqnref{general}, 
under assumptions (A1)--(A3), (H0)--(H3).
Then, strong spectral stability implies 
$L^1\cap L^p\to L^p$ linearized orbital stability, for any $p>1$.
More precisely, 
for initial data $U_0\in L^1\cap L^p$, the solution $U=(u,v)^t(x,t)$
of \eqnref{linearized} satisfies the linear decay bounds
$$
|U(\cdot, t)+\delta(t)\bar U'(\cdot)|_{L^p}\le 
C(1+t)^{-\frac{1}{2}(1-1/p)}(|U_0|_{L^1}+ |U_0|_{L^p}).
\eqnlbl{lindecaybound}
$$
Moreover, provided $U_0\in W^{1,p}$, there hold also the derivative bounds
$$
|\big(v(\cdot, t)+\delta(t)\bar v'(\cdot)\big)_x|_{L^p}\le 
C(1+t)^{-\frac{1}{2}(1-1/p)}(|U_0|_{L^1}+ |U_0|_{W^{1,p}}).
\eqnlbl{dlindecaybound}
$$
\endproclaim

{\bf Proof.}
For \eqnref{lindecaybound},
it is equivalent to show that,
for initial data $G_0\in L^1\cap L^p$, the solution $G(x,t)$
of \eqnref{linsymmetrizable} satisfies 
$$
|G(\cdot, t)-\varphi(\cdot, t)|_{L^p}\le 
C(1+t)^{-\frac{1}{2}(1-1/p)}(|G_0|_{L^1}+ |G_0|_{W^{1,p}}).
\eqnlbl{Gdecaybound}
$$
But, this follows immediately from \eqnref{convolution} 
and bounds \eqnref{tGbounds} and \eqnref{Hbounds}, with $q=p$.
Likewise, for \eqnref{dlindecaybound},
it is equivalent to show that,
for initial data $G_0\in L^1\cap W^{1,p}$, the solution $G(x,t)$
of \eqnref{linsymmetrizable} satisfies 
$$
|\left(\Pi_2(A^0)^{-1}\big(G(\cdot, t)-\varphi(\cdot, t)\big)\right)_x|_{L^p}\le 
C(1+t)^{-\frac{1}{2}(1-1/p)}(|G_0|_{L^1}+ |G_0|_{W^{1,p}}).
\eqnlbl{Gdecaybound}
$$
This follows, similarly, from the derivative analog
$$
\left(\Pi_2 (A^0)^{-1}
\big(G(x,t)-\varphi(x,t) \big) \right)_x=
\int_{-\infty}^{+\infty}
\left(\Pi_2 (A^0)^{-1}(H + \tG)\right)_x(x,t;0)G_0(y)\, dy
\eqnlbl{dconvolution}
$$
of \eqnref{convolution},
together with \eqnref{good1} and bounds \eqnref{tGbounds}, \eqnref{tGxbounds},
with $q=p$.
\myqed
\bigskip

\newsection {\bf Section 5. Nonlinear stability.}
\sectionnumber=5
\theoremnumber=0
\equationnumber=0
\smallskip
\TagsOnLeft
\bigskip

Finally, we establish our main result
of nonlinear orbital stability with respect to perturbations
$ U_0\in L^1\cap H^3 $
of weak viscous profiles (necessarily Lax type, by Proposition \thmref{freist1})
of dissipative, symmetric hyperbolic--parabolic systems
of type \eqnref{general},
of strength $\varepsilon>0$ sufficiently small with respect to the parameters of
the system in question.
We follow the basic iteration scheme of [MZ.1,Z.2];
for precursors of this scheme,  see [Go.2,K.1--2,LZ.1--2,ZH,HZ.1--2]. 

Define the nonlinear perturbation
$$
U(x,t):= \tilde U (x+\delta(t),t)-
\bar U(x),
\eqnlbl{pert}
$$
where $\delta(t)$ (estimating shock location) is to be determined
later; for definiteness, fix $\delta(0)=0$.  
Substituting \eqnref{pert} into \eqnref{general}, we obtain 
$$
G(\tilde U)_t + F(\tilde U)_x - (B(\tilde U)\tilde U_x)_x=
\dot \delta G(\tilde U)_x,
\eqnlbl{tildeU}
$$
and thereby the basic perturbation equation
$$
\left(G(\tilde U)-G(\bar U)\right)_t 
+\left(F(\tilde U)-F(\bar U)\right)_x 
-\left(B(\tilde U)\tilde U_x-B(\bar U)\bar U_x\right)_x 
=\dot \delta(t)G(\tilde U)_x,
\eqnlbl{basicpert}
$$
where $\tilde U$ now denotes $\tilde U(x+\delta(t),t)$ and
$\bar U$ denotes $\bar U(x)$.

\medskip
{\bf 5.1. Energy estimates.}
We begin by establishing the following basic energy estimate, by which we
will eventually close our nonlinear iteration argument:

\proclaim{Proposition \thmlbl{energy}}
Let $U_0\in H^3$, and suppose that, for $0\le t\le T$,
both the supremum of $ |\dot\delta|$ and the $H^2\cap W^{2,\infty}$ 
norm of the solution 
$U=(u,v)^t$ of \eqnref{pert}--\eqnref{basicpert} remain
bounded by a sufficiently
small constant $\zeta>0$, 
for a small-amplitude shock profile as described in Proposition
\thmref{freist1} of a system \eqnref{general} satisfying
(A1)--(A3), (H0)--(H3),
with $\varepsilon :=|\tilde U_+ - \tilde U_-|$ sufficiently small.
Then, there hold the bounds:
$$
|U|_{H^3}(t), \, \int_0^t \big(|U_x|_{H^2}+|v_x|_{H^3}^2\big)(s) ds \le 
C\int_0^t \left(|U|_{L^\infty}(|U|_{L^\infty}+ |U|_{L^2}^2) 
+ \dot \delta^2 \right)(s)ds,
\eqnlbl{ebounds}
$$
for all $0\le t\le T$.
\endproclaim

{\bf Remark \thmlbl{nearoptimal}.}
Note, with the expected decay rates $|U|_{L^\infty}$,
$|\dot \delta(t)|\sim C(1+t)^{-1/2}$,
and $|U|_{L^2}\sim C(1+t^{-1/4})$,
that the righthand side of \eqnref{ebounds}
becomes order $\log (1+t)$, very nearly recovering
the order one bound available in the constant-coefficient
case [Kaw].
\medskip

\proclaim{Lemma \thmlbl{skew} ([SK])}
Assuming (A1), condition (A2)
is equivalent to either of:
\medskip
(K1) \quad 
There exists a smooth skew-symmetric matrix $K(u)$ such that
$$
\text{\rm Re }\left( K(A^0)^{-1}A + B \right)(u) \ge \theta>0.
\eqnlbl{skew}
$$
$A^0$, $A$, $B$ as in \eqnref{AB}.

\medskip  
(K2)  \quad
For some $\theta>0$, there holds
$$
\text{\rm Re }\sigma(-i\xi (A^0)^{-1}A(u)-|\xi|^2 (A^0)^{-1}B(u))
\le -\theta |\xi|^2/(1+|\xi|^2),
\eqnlbl{symbol}
$$
for all $\xi\in \BbbR$.
\endproclaim

{\bf Proof.}
These and other useful equivalent formulations are established in [SK].  \myqed

\proclaim{Lemma \thmlbl{pbounds}}
Let (A1)--(A3) and (H1)--(H3) hold,
and let $\bar U(x-st)$ be a viscous shock solution such that
the profile $\{\bar U(z)\}$ lies entirely within a sufficently 
small neighborhood $\CalV\subset \CalU$ of $\tilde U_*$, and the speed 
$s$ lies within a sufficiently small neighborhood of $a_p(\tilde U_*)$:
i.e., a profile as described in Proposition \thmref{freist1}.
Then, letting $\varepsilon :=|\tilde U_+ - \tilde U_-|$ denote shock strength, 
we have 
for $q=0,\dots 4$ the uniform bounds:
$$
\aligned
|\partial_x^q \bar U(x)| &\le  C\varepsilon^{q+1}e^{-\theta\varepsilon|x|},\\
|\partial_x^q \bar U|_{L^p} &\le  C\varepsilon^{q+1-1/p},\\
\endaligned
\eqnlbl{weakbounds}
$$
for some $C$, $\theta>0$.
\endproclaim

{\bf Proof.}
Though the bounds \eqnref{weakbounds} are not explicitly stated
in [Fre.2], they follow immediately from the detailed description 
of center manifold dynamics obtained in the proof, exactly as
in the strictly parabolic case [MP].  
\myqed

{\bf Proof of Proposition \thmref{energy}}.
Writing \eqnref{pert} in quasilinear form
$$
\left(\tilde A^0\tilde U_t-\bar A^0 \bar U_t \right)
+\left(\tilde A\tilde U_x-\bar A\bar U_x\right)
-\left(\tilde B\tilde U_x-\bar B\bar U_x\right)_x 
=\dot \delta(t)\tilde A^0\tilde U_x,
\eqnlbl{basicpert}
$$
where 
$$
\tilde A^0:=A^0(\tilde U), \bar A^0:=A^0(\bar U);
\quad \tilde A:=A(\tilde U), \bar A:=A(\bar U);
\quad \tilde B:=B(\tilde U), \bar B:=B(\bar U),
\eqnlbl{tildebar}
$$
using the quadratic Leibnitz relation
$$
A_2 U_2- A_1 U_1= A_2(U_2-U_1) + (A_2-A_1)U_1,
\eqnlbl{nLeib}
$$
and recalling the block structure assumption (A3),
we obtain the alternative perturbation equation: 
$$
\tilde A^0 U_t +\tilde A U_x -(\tilde B U_x )_x =
M_1(U)\bar U_x 
+ (M_2(U) \bar U_x)_x
+\dot \delta(t)\tilde A^0 U_x
+\dot \delta(t)\tilde A^0\bar U_x,
\eqnlbl{Leibnitz}
$$
where
$$
M_1(U):=\tilde A-\bar A
=
\left( \int_0^1 dA(\bar U + \theta U )\, d\theta  \right) U,
\eqnlbl{M1}
$$
and
$$
\aligned
M_2(U)&:= \tilde B-\bar B\\
&=
\pmatrix
0&0\\
0& 
( \int_0^1 db(\bar U + \theta U )\, d\theta  ) U
\endpmatrix,
\endaligned
\eqnlbl{M2}
$$
with
$$
\tilde b:=b(\tilde U), \, \bar b:=b(\bar U).
\eqnlbl{btilde}
$$

We now carry out a series of successively higher order energy estimates of the
the type formalized by Kawashima [Ka.1].  
The origin of this approach goes back
to [K,MN] in the context of gas dynamics; see, e.g., [HoZ] for further
discussion/references.

Let $\tilde K$ denote the skew-symmetric matrix 
described in Lemma \thmref{skew}
associated with $\tilde A^0$, $\tilde A$, $\tilde B$.
Then, regarding $\tilde A^0$, $\tilde K$, we have the bounds
$$
\aligned
\tilde A^0_x&=
dA^0(\tilde U) \tilde U_x, \quad
\tilde K_x=
dK(\tilde U) \tilde U_x, \quad
\tilde A_x=
dA(\tilde U) \tilde U_x, \quad
\tilde B_x=
dB(\tilde U) \tilde U_x, 
\\
\tilde A^0_t&=
dA^0(\tilde U) \tilde U_t , \quad
\tilde K_t=
dK(\tilde U) \tilde U_t, \quad
\tilde A_t=
dA(\tilde U) \tilde U_t, \quad
\tilde B_x=
dB(\tilde U) \tilde U_t, \\
\endaligned
\eqnlbl{AKbounds}
$$
and (from defining equations \eqnref{pert}--\eqnref{tildeU}): 
$$
|\tilde U_x|=|U_x + \bar U_x|\le |U_x|+ |\bar U_x|
\eqnlbl{xbound}
$$
and
$$
\aligned
|\tilde U_t|&
\le C( |\tilde U_x| + |\tilde v_{xx}| + |\dot \delta||\tilde U_x|)\\
&
\le C( |U_x| + |\bar U_x| + |v_{xx}|+ |\bar v_{xx}| + |\dot \delta||U_x|
+|\dot\delta| |\bar U_x|)\\
&
\le C( |U_x| + |\bar U_x| + |v_{xx}|+ |\bar v_{xx}| ).\\
\endaligned
\eqnlbl{tbound}
$$
Thus, in particular
$$
|\dot\delta|, \,
|\tilde A^0_x|,\,
|\tilde K_x|,\,
|\tilde A_x|,\,
|\tilde B_x| ,\,
|\tilde A^0_t|,\,
|\tilde K_t|,\,
|\tilde A_t|,\,
|\tilde B_t| 
\le C(\zeta+\varepsilon).
\eqnlbl{smallness}
$$

\medskip

{\it $H^1$ estimate.}  
We first perform a standard, ``Friedrichs-type'' 
estimate for symmetrizable hyperbolic systems.
Taking the $L^2$ inner product of $U$ against \eqnref{Leibnitz}, 
we obtain after rearrangement/integration by parts, and several
applications of Young's inequality, the energy estimate
$$
\aligned
\frac{1}{2}\langle  U_{},\tilde A^0 U_{}\rangle_t 
&=
\langle U_{},\tilde A^0 U_{t}\rangle
+\frac{1}{2}\langle  U_{},\tilde A^0_t U_{}\rangle
\\ 
&=
-\langle U,\tilde A U_x\rangle
  +
\langle U,(\tilde B U_x )_x \rangle
+\langle U, M_1(U)\bar U_x \rangle 
+ \langle U, (M_2(U) \bar U_x)_x\rangle\\
&\qquad+\dot \delta(t)\langle U,\tilde A^0 U_x\rangle
+\dot \delta(t)\langle U,\tilde A^0\bar U_x\rangle
+\frac{1}{2}\langle  U_{},\tilde A^0_t U_{}\rangle\\
&=
\frac{1}{2}\langle U,\tilde A_x U\rangle
  -
\langle U_x,\tilde B U_x  \rangle
+\langle U, M_1(U)\bar U_x \rangle 
- \langle U_x, M_2(U) \bar U_x\rangle\\
&\qquad-\frac{1}{2}\dot \delta(t)\langle U,\tilde A^0_x U\rangle
+\dot \delta(t)\langle U,\tilde A^0\bar U_x\rangle
+\frac{1}{2}\langle  U_{},\tilde A^0_t U_{}\rangle\\
&\le - 
\langle U_{x},\tilde BU_{x}\rangle \\
&\qquad
+ C\int \left( (|U_x|+|\bar U_x|+ |v_{xx}|)|U|^2 
+ |v_x||U||U_x| + |\dot \delta||U||\bar U_x| \right) \\
&\le - 
\langle U_{x},\tilde BU_{x}\rangle \\
&\qquad
+ C\int \left( (|U_x|^2+|U|^2+ |v_x|^2+ |v_{xx}|^2)(|U| + |\bar U_x|) 
+ |\dot \delta|^2|\bar U_x| \right) \\
&\le - 
\langle U_{x},\tilde BU_{x}\rangle \\
&\qquad
+ C\left(|U|_{L^\infty}(|U|_{L^\infty}+ |U|_{L^2}^2) 
+ |\dot \delta|^2 \right)
+ C(\varepsilon +\zeta )\left( |U_x|_{L^2}^2 + |v_{xx}|_{L^2}^2\right),
\endaligned
\eqnlbl{sym}
$$
where $\varepsilon$, $\zeta>0$ is as in the statement of the Proposition.
Here, we have freely used the weak shock assumption and 
consequent bounds \eqnref{weakbounds}, as well as \eqnref{AKbounds}.
(Note: we have also used in a crucial way the block-diagonal
form of $M_2$ in estimating
$|\langle U_x, M_2(U) \bar U_x\rangle|\le 
C\int |v_x||U||\bar U_x|$ in the first inequality).

Likewise, differentiating \eqnref{Leibnitz}, taking the $L^2$
inner product of $U_x$ against the resulting equation,
and substituting the result into
$$
\frac{1}{2}\langle U_{x},\tilde A^0 U_{x}\rangle_t 
=
\langle U_{x},(\tilde A^0 U_{t})_x\rangle
-\langle U_{x},\tilde A^0_x U_{t}\rangle
+\frac{1}{2}\langle  U_{x},\tilde A^0_t U_{x}\rangle,
\eqnlbl{one}
$$
we obtain after an integration by parts: 
$$
\aligned
\frac{1}{2}\langle U_{x},\tilde A^0 U_{x}\rangle_t 
&=
-\langle U_x,(\tilde A U_x)\rangle
  +
\langle U,(\tilde B U_x )_x \rangle
+\langle U, M_1(U)\bar U_x \rangle 
+ \langle U, (M_2(U) \bar U_x)_x\rangle\\
&\qquad+\dot \delta(t)\langle U,\tilde A^0 U_x\rangle
+\dot \delta(t)\langle U,\tilde A^0\bar U_x\rangle
+\frac{1}{2}\langle  U_{},\tilde A^0_t U_{}\rangle\\
&=
-\frac{1}{2}\langle U_x,\tilde A_x U_x\rangle
- \langle U_{xx},\tilde B U_{xx}  \rangle
- \langle U_{xx},\tilde B_x U_{x}  \rangle
+\langle U_x, (M_1(U)\bar U_x)_x \rangle \\
&\qquad - \langle U_{xx}, (M_2(U) \bar U_x)_x\rangle
+\frac{1}{2}\dot \delta(t)\langle U_x,\tilde A^0_x U_x\rangle
+\dot \delta(t)\langle U_x,\tilde A^0\bar U_{xx}\rangle\\
&\qquad +\dot \delta(t)\langle U_x,\tilde A^0_x\bar U_x\rangle
-\langle U_{x},\tilde A^0_x U_{t}\rangle
+\frac{1}{2}\langle  U_{x},\tilde A^0_t U_{x}\rangle,\\
\endaligned
\eqnlbl{nextsym0}
$$
which by \eqnref{smallness}, plus various applications of Young's
inequality yields the next-order energy estimate:
$$
\aligned
\frac{1}{2}\langle U_{x},\tilde A^0 U_{x}\rangle_t 
&\le - 
\langle U_{xx},\tilde BU_{xx}\rangle \\
&\qquad
+ C\int\left(
 (|U|^2+|\dot\delta|^2)(|\bar U_{xx}|+ |\bar U_x|)
+ (\varepsilon +\zeta )( |U_x|^2 + |U_x||v_{xx}| ) \right),\\
&\le - 
\langle U_{xx},\tilde BU_{xx}\rangle \\
&\qquad
+ C\left((|U|_{L^\infty}^2+ 
 |\dot \delta|^2 \right)
+ C(\varepsilon +\zeta )\left( |U_x|_{L^2}^2 + |v_{xx}|_{L^2}^2\right)\\
&\le - \frac{1}{2}
\langle U_{xx},\tilde BU_{xx}\rangle \\
&\qquad
+ C\left((|U|_{L^\infty}^2+ 
 |\dot \delta|^2 \right)
+ C(\varepsilon +\zeta )|U_x|_{L^2}^2.
\endaligned
\eqnlbl{nextsym}
$$

Next, we perform a nonstandard, ``Kawashima-type'' derivative estimate.
Taking the $L^2$ inner product of $U_{x}$
against $\tilde K(\tilde A^0)^{-1} $ times \eqnref{Leibnitz},
and noting that (integrating by parts, and using skew-symmetry of $\tilde K$)
$$
\aligned
\frac{1}{2}\langle U_{x},\tilde KU_{} \rangle_t
&=
 \frac{1}{2}\langle U_{x},\tilde KU_{t} \rangle+
\frac{1}{2}\langle U_{xt},\tilde KU_{} \rangle
+\frac{1}{2}\langle U_{x},\tilde K_tU_{} \rangle\\
&=
\frac{1}{2}\langle U_{x},\tilde KU_{t} \rangle
-\frac{1}{2}\langle U_{t},\tilde KU_{x} \rangle\\
&\quad
-\frac{1}{2}\langle  U_{t},\tilde K_x U_{x} \rangle
+\frac{1}{2}\langle U_{x},\tilde K_tU_{} \rangle\\
&=
\langle U_{x},\tilde KU_{t} \rangle
+\frac{1}{2}\langle U_{},\tilde K_x U_{t} \rangle
+\frac{1}{2}\langle U_{x},\tilde K_tU_{} \rangle ,\\
\endaligned
\eqnlbl{parts}
$$
we obtain by calculations similar to the above
the auxiliary energy estimate:
$$
\aligned
\frac{1}{2}\langle U_{x},\tilde KU_{} \rangle_t &\le
-\langle  U_{x}, \tilde K(\tilde A^0)^{-1}\tilde A U_{x}\rangle\\
&\qquad
 + C(\zeta+\varepsilon) |U_x|_{L^2}^2 + 
C\zeta^{-1}|v_{xx}|_{L^2}^2 + C(|U|_{L^\infty}^2+ |\dot\delta(t)|^2).
\endaligned
\eqnlbl{est2}
$$
Adding \eqnref{sym}, \eqnref{est2}, and 
\eqnref{nextsym} times a suitably large constant $M>0$, 
and recalling \eqnref{skew},
we obtain, finally:
$$
\aligned
\frac{1}{2}
\Big(
\langle U, \tilde A^0 U_{}\rangle
+ &\langle U_{x},\tilde KU_{} \rangle
+M\langle U_x, \tilde A^0 U_{x}\rangle \Big)_t \\
&\le
-\theta (|U_{x}|_{L^2}^2
+|v_{xx}|_{L^2}^2)
+C \left(|U|_{L^\infty}(|U|_{L^\infty}+ |U|_{L^2}^2) 
+ |\dot \delta|^2 \right),\\
\endaligned
\eqnlbl{estfinal}
$$
$\theta>0$,
for any $\zeta$, $\varepsilon $ sufficiently small, and $M$, $C>0$ sufficiently large.

\medskip
{\it Higher order estimates.}
Performing the same procedure on the once- and twice-differentiated
versions of equation \eqnref{Leibnitz}, we obtain, 
likewise, the $H^q$ estimates,
$q=2, \, 3$, of:
$$
\aligned
\frac{1}{2}
\Big(
\langle  \partial_x^{q-1}U_{},&\tilde A^0\partial_x^{q-1}U_{}\rangle
+ \langle \partial_x^{q-1}U_{x},\tilde K\partial_x^{q-1}U_{} \rangle
+M\langle  \partial_x^{q}U,\tilde A^0\partial_x^{q}U\rangle \Big)_t \\
&\le
-\theta (|\partial_x^q U|_{L^2}^2
+|\partial_x^{q+1}v|_{L^2}^2)\\
&\qquad +(\varepsilon + \zeta) |U_{x}|_{H^{q-2}}^2
+C \left(|U|_{L^\infty}(|U|_{L^\infty}+ |U|_{L^2}^2) 
+ |\dot \delta|^2 \right).\\
\endaligned
\eqnlbl{esthigher}
$$
We omit the calculations, which are entirely similar to those
carried out already.
\medskip

{\it Final estimate.}
Summing our $H^q$ estimates from $q=1$ to $3$,
and telescoping the sum of the righthand sides,
we thus obtain, for $\varepsilon$, $\zeta$ sufficiently small:
$$
\aligned
\sum_{q=1}^3 
\big( \langle &\tilde A^0 \partial_x^{q-1}U,\partial_x^{q-1}U\rangle
+ \frac{1}{2}\langle \partial_x^{q}U,\tilde K\partial_x^{q-1}U
\rangle
+M\langle \tilde A^0 \partial_x^{q}U,\partial_x^{q}U\rangle \big)_t \\
&\le
-\theta (|U_x|_{H^2}+|v_x|_{H^3}^2)
+C \left(|U|_{L^\infty}(|U|_{L^\infty}+ |U|_{L^2}^2) 
+ |\dot \delta|^2 \right),\\
\endaligned
\eqnlbl{estsummed}
$$
or,
integrating from $0$ to $t$: 
$$
\aligned
\sum_{q=1}^3 
\big( \langle \tilde A^0 \partial_x^{q-1}U,\partial_x^{q-1}U\rangle
+ &\frac{1}{2}\langle \partial_x^{q}U,\tilde K\partial_x^{q-1}U
\rangle
+M\langle \tilde A^0 \partial_x^{q}U,\partial_x^{q}U\rangle \big)|^t_0 \\
&\le
-\int_0^t \theta (|U_x|_{H^2}+|v_x|_{H^3}^2)(s)ds\\
&\quad+
C\int_0^t \left(|U|_{L^\infty}(|U|_{L^\infty}+ |U|_{L^2}^2) 
+ |\dot \delta|^2 \right)(s)ds.
\endaligned
\eqnlbl{intcalc0}
$$

Noting that, for $M$ sufficiently large, we have for each $q$,
by Young's inequality, and positive definiteness of $\tilde A^0$:
$$
\big(
\langle  \partial_x^{q-1}U,\tilde A^0 \partial_x^{q-1}U\rangle
+ \frac{1}{2}\langle \partial_x^{q}U,\tilde K\partial_x^{q-1}U
\rangle
+M\langle  \partial_x^{q}U,\tilde A^0 \partial_x^{q}U\rangle \big)(t)\ge
\theta(|\partial_x^{q-1} U|_{L^2}^2
+|\partial_x^{q} U|_{L^2}^2),
\eqnlbl{absorb}
$$
for some $\theta>0$, we may rearrange \eqnref{intcalc0} to obtain
our ultimate goal:
$$
\aligned
|U|_{H^3}^2(t)
+&\int_0^t \theta (|U_x|_{H^2}+|v_x|_{H^3}^2)(s)ds\\
&\le  C|U|_{H^3}^2(0)+
C\int_0^t \left(|U|_{L^\infty}(|U|_{L^\infty}+ |U|_{L^2}^2) 
+ |\dot \delta|^2 \right)(s)ds,
\endaligned
\eqnlbl{goal}
$$
from which the result immediately follows.
\myqed

\medskip

{\bf 5.2. Nonlinear iteration.}
We now carry out the nonlinear iteration, following [Z.2,MZ.1].
For this stage of the argument, it will be convenient to work
again with the conservative variable 
$$
G:=G(\tilde U)-G(\bar U),
\eqnlbl{Gpert}
$$ 
writing \eqnref{pert} in the more standard form:
$$
G_t-LG=N(G,G_x)_x 
+\dot \delta (t)(\BG_x + G_x),
\eqnlbl{29}
$$
$\bar G:= G(\bar U)$, where
$$
\aligned
N(G,G_x)&=\Cal{O}(|G|^2 + |G||v_x|),\\
N(G,G_x)_x&=\Cal{O}(|G|^2 + |G||v_x|
+ |G_x||v_x| + |G||v_{xx}|),\\
\endaligned
\eqnlbl{29.1}
$$
so long as $|G|$, $|G_x|$ remain bounded.
Here, $v$ denotes the second coordinate of the alternative 
perturbation variable $U=(u,v)^t$ defined in \eqnref{pert}.

By Duhamel's principle, and the fact that
$$
\int^\infty_{-\infty}\CalG(x,t;y)\BG_x(y)dy=e^{Lt}\BG_x(x)=\BG_x(x),
\eqnlbl{stationary}
$$
we have, formally,
$$
\aligned
& G(x,t)=\int^\infty_{-\infty}\CalG(x,t;y)G_0(y)dy\\
&+\int^t_0 \int^\infty_{-\infty} \CalG_y(x,t-s;y)
(N(G,G_x)+
\dot \delta G ) (y,s)dyds\\
&+ \delta (t)\BG_x.\\
\endaligned
\eqnlbl{30}
$$

Defining, by analogy with the linear case,
the {\it nonlinear instantaneous projection}:
$$
\aligned
\varphi(x,t)
&:= -\delta(t)\BG_x\\
&:= \int^\infty_{-\infty}{E}(x,t;y)G_0(y) dy\\
&-\int^t_0 \int^\infty_{-\infty}E_y(x,t-s;y)(N(G,G_x)+
\dot \delta G)(y,s)dy,\\
\endaligned
\eqnlbl{proj}
$$
or equivalently, the {\it instantaneous shock location}:
$$
\aligned
\delta (t)
&=-\int^\infty_{-\infty}e(y,t) G_0(y)dy\\
&+\int^t_0\int^{+\infty}_{-\infty} e_{y}(y,t-s)(N(G,G_x)+
\dot \delta G)(y,s) dy ds,\\
\endaligned
\eqnlbl{31}
$$
where $E$, $e$ are defined as in \eqnref{E}, \eqnref{e},
and recalling \eqnref{tG}, we thus
obtain the {\it reduced equations}:
$$
\aligned
G(x,t)
&=\int^\infty_{-\infty} (H+\tG)(x,t;y)G_0(y)dy\\
&+\int^t_0
\int^\infty_{-\infty}H(x,t-s;y)(N(G,G_x)_x
+\dot \delta G_x)(y,s)dy \, ds\\
&-\int^t_0
\int^\infty_{-\infty}\tG_y(x,t-s;y)(N(G,G_x)+
\dot \delta G)(y,s) dy \, ds,\\
\endaligned
\eqnlbl{32}
$$
and, differentiating \eqnref{31} with respect to $t$,
$$
\aligned
\dot \delta (t)
&=-\int^\infty_{-\infty}e_t(y,t) G_0(y)dy\\
&+\int^t_0\int^{+\infty}_{-\infty} e_{yt}(y,t-s)(N(G,G_x)+
\dot \delta G)(y,s) dy ds.\\
\endaligned
\eqnlbl{33}
$$

{\it Note:}  
In deriving \eqnref{33}, we have used the fact 
that $e_y (y,s)\rightharpoondown 0$ as $s \to 0$, as the 
difference of approaching heat kernels, in evaluating the 
boundary term
$$
\int^{+\infty}_{-\infty} e_y (y,0)(N(G,G_x)+
\dot \delta G)(y,t)dy=0.
\eqnlbl{35}
$$
(Indeed, $|e_y(\cdot, s)|_{L^1}\to 0$, see Remark 2.6, below).
\medskip

The defining relation $ \delta (t)\bar u_x:= -\varphi$ in \eqnref{proj}
can be motivated heuristically by
$$
\eqalign{
\tilde G(x,t)-\varphi(x,t) \sim G&= \pmatrix u\\v\endpmatrix (x+\delta(t),t) 
- \pmatrix \bar u\\ \bar v\endpmatrix(x) \cr
&\sim \tilde G(x,t) + \delta (t)\bar G_x(x),}
$$
where $\tilde G$ denotes the solution of the linearized
perturbation equations, and $\bar G$ the background profile.  
Alternatively, it can be thought
of as the requirement that the instantaneous projection 
of the shifted (nonlinear) perturbation variable $G$ be zero, [HZ.1--2].

\proclaim{Lemma \thmlbl{3} [Z.2]}  The kernel ${e}$ satisfies
$$
|{e}_y (\cdot, t)|_{L^p},  |{e}_t(\cdot, t)|_{L^p} 
\le C t^{-\frac{1}{2}(1-1/p)},
\eqnlbl{36}
$$
$$
|{e}_{ty}(\cdot, t)|_{L^p} 
\le C t^{-\frac{1}{2}(1-1/p)-1/2},
\eqnlbl{37}
$$
for all $t>0$.  Moreover, for $y\le 0$ we have the pointwise bounds
$$
|{e}_y (y,t)|, |{e}_t (y,t)| \le Ct^{-\frac{1}{2}} e^{-\frac{(y+a_-t)^2}{Mt}},
\eqnlbl{38}
$$
$$
|{e}_{ty} (y,t)| \le C t^{-1} e^{-\frac{(y+a_-t)^2}{Mt}}, 
\eqnlbl{39}
$$
for $M>0$ sufficiently large (i.e. $>4b_\pm$),
and symmetrically for $y\ge 0$.
\endproclaim

{\bf Proof.}  
For definiteness, take $y \le 0$. 
Then, \eqnref{e} gives
$$
{e}_y(y,t)=
\left(\frac{1}{u_+-u_-}\right)
\left(K(y+a_-t,t)-K(y-a_-t,t)\right),
\eqnlbl{41}
$$
$$
{e}_t(y,t)=
\left(\frac{1}{u_+-u_-}\right)
\left( (K+K_y)(y+a_-t,t)-(K+K_y)(y-a_-t,t)\right),
\eqnlbl{41}
$$
$$
{e}_{ty}(y,t)=
\left(\frac{1}{u_+-u_-}\right)
\left( (K_y+K_{yy})(y+a_-t,t)-(K_y+K_{yy})(y-a_-t,t)\right),
\eqnlbl{42}
$$
where
$$
K(y,t):= \frac{e^{-y^2/4b_-t}}{\sqrt{4\pi b_-t}}
\eqnlbl{43}
$$
denotes an appropriate heat kernel.  The pointwise bounds
\eqnref{38}--\eqnref{39} follow immediately for $t\ge 1$ 
by properties of the heat kernel, in turn yielding \eqnref{36}--\eqnref{37} 
in this case.  The bounds for small time $t\le 1$ follow from estimates
$$
\aligned
|K_y (y+a_-t,t)-K_y (y-a_-t,t)|
&=|\int^{y-a_-t}_{y+a_-t}K_{yy}(z,t)dz| \\
&\le Ct^{-3/2}\int^{y-a_-t}_{y+a_-t} e^{\frac{-z^2}{Mt}} dz\\
&\le Ct^{-1/2}e^{-\frac{(y+a_-t)^2}{Mt}},\\
\endaligned
\eqnlbl{44}
$$
and, similarly,
$$
\aligned
|K_{yy}(y+a_-t,t)-K_{yy}(y-a_-,t)|
&= |\int^{a_-t}_{-a_-t}K_{yyy}(z,t)dz|\\
&\le Ct^{-2}\int^{y-a_-t}_{y+a_-t} e^{\frac{-z^2}{Mt}} dz,\\
&\le Ct^{-1}e^{-\frac{(y+a_-t)^2}{Mt}}.\\
\endaligned
\eqnlbl{45}
$$
The bounds for $|{e}_y|$ are again immediate.
Note that we have taken crucial account of cancellation in
the small time estimates of $e_t$, $e_{ty}$.
\qed

{\bf Remark \thmlbl{2.6}:}
For $t\le 1$, a calculation analogous to that of
\eqnref{44} yields
$
|e_y(y,t)|\le 
Ce^{-\frac{(y+a_-t)^2}{Mt}},
$
and thus $|e(\cdot,s)|_{L^1}\to 0$ as $s\to 0$.

\medskip

With these preparations, we are ready to carry out our analysis:
\medskip

{\bf Proof of Theorem \thmref{nonlin}.} 
Define
$$
\aligned
\zeta(t)
&:= \sup_{0\le s \le t,\, 2\le p\le \infty}
\Big[ \, \big(|U(\cdot, s)|_{L^p}+ |v_x(\cdot, s)|_{L^p}\big) 
(1+s)^{\frac{1}{2}(1-1/p)}\\
&\qquad
+ |\dot \delta (s)|(1+s)^{1/2} 
+ \big( |\delta (s)|
+ |U(\cdot, s)|_{H^2}
+ |U(\cdot, s)|_{W^{2,\infty}}\big)
\, \Big].\\
\endaligned
\eqnlbl{zeta2}
$$
We shall establish:

{\it Claim.} For all $t\ge 0$ for which a solution exists with
$\zeta$ uniformly bounded by some fixed, sufficiently small constant,
there holds
$$
\zeta(t) \leq C_2(\zeta_0 + \zeta(t)^2).
\eqnlbl{claim}
$$
\medskip
{}From this result, it follows by continuous induction that,
provided $\zeta_0 < 1/4C_2$,  
there holds 
$$
\zeta(t) \leq 2C_2\zeta_0
\eqnlbl{bd}
$$
for all $t\geq 0$ such that $\zeta$ remains small.
By standard short-time theory/continuation, we find that
the solution (unique, in this regularity class) in fact remains in 
$H^2$ for all $t\ge 0$,
with bound \eqnref{bd}, at once yielding existence and the claimed 
sharp $L^p$ bounds, $2\le p\le \infty$.
Thus, it remains only to establish the claim above.
\medskip

{\it Proof of Claim.}
We must show that each of the quantities $|U|_{L^p}(1+s)^{\frac{1}{2}(1-1/p)}$,
$|v_x|_{L^p} (1+s)^{\frac{1}{2}(1-1/p)}$,
$|\dot \delta|(1+s)^{1/2}$, $|\delta|$, $|U|_{H^2}$,
and $|U|_{W^{2,\infty}}$ are separately bounded by
$$
C(\zeta_0 + \zeta(t)^2),
\eqnlbl{sep}
$$
for some $C>0$, all $0\le s\le t$, so long as $\zeta$ remains
sufficiently small.

Provided that we can establish the others, the
final two bounds follow easily from the energy estimates of
Proposition \thmref{energy}.
For, by bounds \eqnref{sep} on $|\dot\delta|$ and $|U|_{L^\infty}$,
and the assumption that $|U|_{H^2\cap W^{2,\infty}} \le \zeta$
remains small, we obtain from \eqnref{ebounds} that
$$
|U|_{H^3}(s)\le 
C_3(\zeta_0 + \zeta(t)^2)\log (1+S),
\eqnlbl{logbound}
$$
for $0\le s\le t$.
Interpolating with assumed bound 
$$
|U|_{L^2}(s)\le 
C_3(\zeta_0 + \zeta(t)^2)(1+s)^{-1/4},
$$
we obtain
$$
\aligned
|U|_{H^2}(s)&\le 
C_3(\zeta_0 + \zeta(t)^2)|U|_{H^3}^{2/6}|U|_{L^2}^{1/6}\\
&\le 
C_3(\zeta_0 + \zeta(t)^2)(1+t)^{-1/25},
\endaligned
\eqnlbl{H2bound}
$$
and, by Sobolev estimate,
$$
|U|_{W^{2,\infty}}\le
|U|_{H^2}^{1/2} |U|_{H^3}^{1/2}
\le 
C_3(\zeta_0 + \zeta(t)^2)(1+t)^{-1/51},
\eqnlbl{Wkpbound}
$$
both better than claimed.

Thus, in order to establish the result, 
we have only to establish the remaining bounds,
on $|U|_{L^p}$, $|v_x|_{L^p}$,
$|\dot \delta|$, and $|\delta|$.
These will be carried out using the Green's function
estimates of Lemma \thmref{2.05}.
Accordingly, we first convert the problem to conservative,
$G$ coordinates, via:

{\it Observation \thmlbl{Gequiv}.}
It is sufficient to establish corresponding bounds on 
$|G|_{L^p}$, $|v_{G,x}|_{L^p}$,
$|\dot \delta|$, and $|\delta|$, where
$
v_G(x,t):=A^0(x)^{-1}G(x,t).
$
\medskip

{\it Proof.}
We have
$$
\aligned
U&=G^{-1}(\tilde G) - G^{-1}(\bar G)
= (A^0)^{-1}_{\text{ave}}(x,t)G\\
&:=\left(\int_0^1 (A^0)^{-1}\big(\bar U(x) + \theta U(x,t)\big)d\theta
\right) \, G,
\endaligned
\eqnlbl{roe}
$$
where 
$$
|(A^0)^{-1}_{\text{ave}}(x,t) -A^0(x)^{-1}|
\le 
C|U|
\le
C|G|
\eqnlbl{linfty}
$$
and
$$
|\left((A^0)^{-1}_{\text{ave}}(x,t) -A^0(x)^{-1}\right)_x|
\le 
C(|\bar U_x|+ |U_x|)
\le 
\zeta,
\eqnlbl{Winfty}
$$
whence
$$
|U|_{L^p}\le C|G|_{L^p}
$$
and
$$
|U_x - ((A^0)^{-1}G)_x|_{L^p}\le
C(|G|_{L^p}|G_x|_{L^\infty} + \zeta |G|_{L^p})
\le
C\zeta |G|_{L^p},
$$
from which the result easily follows.
\myqed

By \eqnref{32}--\eqnref{33}, we have
$$
\aligned
|G|_{L^p}(t)&\le
|\int^\infty_{-\infty} (H+\tG)(x,t;y)G_0(y)dy|_{L^p}\\
&+|\int^t_0
\int^\infty_{-\infty}H(x,t-s;y)(N(G,G_x)_x
+\dot \delta G_x)(y,s)dy \, ds |_{L^p}\\
&+|\int^t_0
\int^\infty_{-\infty}\tG_y(x,t-s;y)(N(G,G_x)+
\dot \delta G)(y,s) dy \, ds|_{L^p}\\
&=: I_a + I_b + I_c ,\\
\endaligned
\eqnlbl{reduced}
$$
$$
\aligned
|v_{G,x}|_{L^p}(t)
&:=
|\big((A^0)^{-1}G\big)_x|_{L^p}(t)
\le
|\int^\infty_{-\infty} \left((A^0)^{-1}\tG\right)_x(x,t;y)G_0(y)dy|_{L^p}\\
&+|\int^t_0
\int^\infty_{-\infty}\big(
(A^0)^{-1} \tG_y\big)_x(x,t-s;y)(N(G,G_x)+
\dot \delta G)(y,s) dy \, ds|_{L^p}\\
&=: II_a + II_b ,\\
\endaligned
\eqnlbl{vreduced}
$$
$$
\aligned
|\dot \delta| (t)
&\le |\int^\infty_{-\infty}e_t(y,t) G_0(y)dy|\\
&+|\int^t_0\int^{+\infty}_{-\infty} e_{yt}(y,t-s)(N(G,G_x)+
\dot \delta G)(y,s) dy ds|\\
&=: III_a + III_b,\\
\endaligned
\eqnlbl{deltadot}
$$
and
$$
\aligned
|\delta| (t)
&\le |\int^\infty_{-\infty}e(y,t) G_0(y)dy|\\
&+|\int^t_0\int^{+\infty}_{-\infty} e_{y}(y,t-s)(N(G,G_x)+
\dot \delta G)(y,s) dy ds|\\
&=: IV_a + IV_b.\\
\endaligned
\eqnlbl{delta}
$$

We estimate each term in turn, following the approach of [Z.2,MZ.1].
The linear term $I_a$  satisfies bound
$$
I_a \le  C\zeta_0 (1+t)^{-\frac{1}{2}(1-1/p)},
\eqnlbl{Ia}
$$ 
as already shown in
the proof of Corollary \thmref{suff}.
Likewise, applying the bounds of Lemma \thmref{2.05} together
with \eqnref{29.1} and definition \eqnref{zeta2}, we have
$$
\aligned
I_b&= 
|\int^t_0 \int^\infty_{-\infty}H(x,t-s;y)(N(G,G_x)_x
+\dot \delta G_x)(y,s)dy \, ds |_{L^p}\\
&\le
C\int_0^t e^{-\eta (t-s)}(|G|_{L^\infty}+|v_{G,x}|_{L^\infty}
+ |\dot \delta|)|G|_{W^{2,p}}(s) ds\\
&\le
C\zeta(t)^2
\int_0^t e^{-\eta (t-s)}(1+s)^{-1/2}ds\\
&\le
C\zeta(t)^2 (1+t)^{-1/2},\\
\endaligned
\eqnlbl{Ib}
$$
and (taking $q=2$ in \eqnref{tGybounds})
$$
\aligned
I_c&= 
|\int^t_0 \int^\infty_{-\infty}\tG_y(x,t-s;y)(N(G,G_x)+
\dot \delta G)(y,s) dy \, ds|_{L^p}\\
&\le
C\int_{0}^t e^{-\eta(t-s)}(|G|_{L^\infty}+|v_{G,x}|_{L^\infty}
+ |\dot \delta|)|G|_{L^p}(s) ds\\
&\quad+
C\int_0^{t} (t-s)^{-3/4+1/2p}
(|G|_{L^\infty}+|v_{G,x}|_{L^\infty}
+ |\dot \delta|)|G|_{L^2}(s) ds\\
&\le
C\zeta(t)^2
\int_{0}^t e^{-\eta(t-s)}
(1+s)^{-\frac{1}{2}(1-1/p)-1/2}ds\\
&\quad +
C\zeta(t)^2
\int_0^{t} (t-s)^{-3/4+1/2p}
(1+s)^{-3/4}ds\\
&\le
C\zeta(t)^2 (1+t)^{-\frac{1}{2}(1-1/p)}.\\
\endaligned
\eqnlbl{Ic}
$$
Summing bounds \eqnref{Ia}--\eqnref{Ic}, we obtain \eqnref{sep},
as claimed, for $2\le p\le \infty$.
The desired bounds on $II_a$ and $II_b$ follow by an identical
calculation, once we notice that, in Lemma \thmref{2.05},
$\left((A^0)^{-1}\tG\right)$ and
$\left((A^0)^{-1}\tG\right)_y$
satisfy the same $L^q\to L^p$  bounds as do $(H+\tG)$ and $(H+\tG)_y$,
respectively.

Similarly, applying the bounds of Lemma \thmref{3} together
with definition \eqnref{zeta2}, we find that
$$
\aligned
III_a&=
\big|\int^\infty_{-\infty}e_t(y,t) G_0(y)dy\big|\\
&\le
|e_t(y,t)|_{L^\infty}(t) |G_0|_{L^1}\\
&\le
C\zeta_0 (1+t)^{-1/2}
\endaligned
\eqnlbl{IIIa}
$$
and 
$$
\aligned
III_b&=
|\int^t_0\int^{+\infty}_{-\infty} e_{yt}(y,t-s)(N(G,G_x)+
\dot \delta G)(y,s) dy ds|\\
&\le
\int^t_0
|e_{yt}|_{L^2}(t-s) 
(|G|_{L^\infty}+|v_{G,x}|_{L^\infty} + |\dot \delta|)
|G|_{L^2}(s) ds\\
&\le
C\zeta(t)^2 
\int^t_0
(t-s)^{-3/4} (1+s)^{-3/4} ds\\
&\le
C\zeta(t)^2 (1+t)^{-1/2},
\endaligned
\eqnlbl{IIIb}
$$
while
$$
\aligned
IV_a&=
\big|\int^\infty_{-\infty}e(y,t) G_0(y)dy\big|\\
&\le
|e(y,t)|_{L^\infty}(t) |G_0|_{L^1}\\
&\le
C\zeta_0 
\endaligned
\eqnlbl{IVa}
$$
and 
$$
\aligned
IV_b&=
\big|\int^t_0\int^{+\infty}_{-\infty} e_{y}(y,t-s)
\dot \delta G(y,s) dy ds\big|\\
&\le
\int^t_0
|e_{y}|_{L^2}(t-s) 
(|G|_{L^\infty}+|v_{G,x}|_{L^\infty} + |\dot \delta|)
|G|_{L^2}(s) ds\\
&\le
C\zeta(t)^2 
\int^t_0
(t-s)^{-1/4} (1+s)^{-3/4} ds\\
&\le
C\zeta(t)^2. 
\endaligned
\eqnlbl{IVb}
$$
Summing \eqnref{IIIa}--\eqnref{IIIb} and
\eqnref{IVa}--\eqnref{IVb},
we obtain \eqnref{sep} as claimed.

This completes the proof of the claim, and the result.
\myqed

\bigskip

\Refs
\medskip\noindent
\medskip\noindent
[CS] C. Conley and J. Smoller,
{\it On the structure of magnetohydrodynamic shock waves,}
Comm. Pure Appl.  Math. 28 (1974), 
367--375.
\medskip\noindent
[Er] J. J. Erpenbeck,
{\it Stability of step shocks,} Phys. Fluids 5 (1962) no. 10, 1181--1187.
\medskip\noindent
[Fre.1] H. Freist\"uhler,
{\it Small amplitude intermediate magnetohydrodynamic shock waves,}
Phys. Scripta T74 (1998) 26--29.
\medskip\noindent
[Fre.2] H. Freist\"uhler,
{\it Profiles for small Laxian shock waves in Kawashima 
type systems of conservation laws,} Preprint (2000).
\medskip\noindent
[FreZ] H. Freist\"uhler and K. Zumbrun,
{\it Examples of unstable viscous shock waves,}
unpublished note, Institut f\"ur Mathematik, RWTH Aachen, February 1998.
\medskip\noindent
[Fri.1] C. Fries, 
{\it Stability of viscous shock waves 
associated with non-convex modes,}
Arch. Ration. Mech. Anal. 152 (2000), no. 2, 141--186.
\medskip\noindent
[Fri.2] C. Fries, 
{\it Nonlinear asymptotic stability of general 
small-amplitude viscous Laxian shock waves,}
J. Differential Equations 146 (1998), no. 1, 185--202.
\medskip\noindent
[GZ] R. Gardner and K. Zumbrun,
{\it The Gap Lemma and geometric criteria for instability
of viscous shock profiles}, 
Comm. Pure Appl.  Math. 51 (1998), no. 7, 797--855. 
\medskip\noindent
[Gi] D. Gilbarg, 
{\it The existence and limit behavior of 
the one-dimensional shock layer,} Amer. J. Math. 73, (1951). 256--274. 
\medskip\noindent
[Go.1] J. Goodman, {\it Nonlinear asymptotic stability of viscous
shock profiles for conservation laws,}
Arch. Rational Mech. Anal. 95 (1986), no. 4, 325--344.
\medskip\noindent
[Go.2] J. Goodman, 
{\it Remarks on the stability of viscous shock waves}, in:
Viscous profiles and numerical methods for shock waves 
(Raleigh, NC, 1990), 66--72, SIAM, Philadelphia, PA, (1991).
\medskip\noindent
[HoZ] D. Hoff and K. Zumbrun, 
{\it Multi-dimensional diffusion waves for the Navier-Stokes equations of
compressible flow,} Indiana Univ. Math. J. 44 (1995), no. 2, 603--676.
\medskip\noindent
[HZ.1] P. Howard and K. Zumbrun, 
{\it A tracking mechanism for one-dimensional 
viscous shock waves,} preprint (1999).
\medskip\noindent
[HZ.2] P. Howard and K. Zumbrun,
{\it Pointwise estimates for dispersive-diffusive shock waves,}
to appear, Arch. Rational Mech. Anal. 
\medskip\noindent
[HuZ] J. Humpherys and K. Zumbrun,
{\it Spectral stability of small amplitude
shock profiles for dissipative symmetric hyperbolic--parabolic systems,}
to appear ZAMP (2001).
\medskip\noindent
[Ka] Ya. Kanel,
{\it On a model system of equations of of one-dimensional
gas motion,} Diff. Eqns. 4 (1968) 374--380.
\medskip\noindent
[K.1] T. Kapitula,
{\it Stability of weak shocks in $\lambda$--$\omega$ systems},
Indiana Univ. Math. J. 40 (1991), no. 4, 1193--12.
\medskip\noindent
[K.2] T. Kapitula,
{\it On the stability of travelling waves in weighted $L^\infty$ spaces},
J. Diff. Eqs. 112 (1994), no. 1, 179--215.
\medskip\noindent
[Kaw] S. Kawashima,
{\it Systems of a hyperbolic--parabolic composite type,
 with applications to the equations of magnetohydrodynamics},
thesis, Kyoto University (1983).
\medskip\noindent
[KM] S. Kawashima and A. Matsumura,
{\it Asymptotic stability of traveling wave solutions of systems 
for one-dimensional gas motion,} Comm.
Math. Phys. 101 (1985), no. 1, 97--127.
\medskip\noindent
[KMN] S. Kawashima, A. Matsumura, and K. Nishihara,
{\it Asymptotic behavior of solutions for the equations of a 
viscous heat-conductive gas,}
Proc.  Japan Acad. Ser. A Math. Sci. 62 (1986), no. 7, 249--252. 
\medskip\noindent
[La] P.D. Lax,
{\it Hyperbolic systems of conservation laws and the mathematical theory of shock waves},
Conference Board of the Mathematical Sciences Regional Conference
Series in Applied Mathematics, No. 11. Society for Industrial and Applied Mathematics,
Philadelphia, Pa., 1973. v+48 pp.
\medskip\noindent
[L.1] T.-P. Liu,
{\it Pointwise convergence to shock waves for viscous conservation laws},
Comm. Pure Appl. Math. 50 (1997), no. 11, 1113--1182.
\medskip\noindent
[L.2] T.-P. Liu, 
{\it Shock waves for compressible Navier--Stokes equations are stable,}
Comm. Pure Appl. Math. 39 (1986), no. 5, 565--594.
\medskip\noindent
[L.3] T.-P. Liu, 
{\it The Riemann problem for general systems of conservation laws,}
J. Diff. Eqs. 18 (1975) 218--234.
\medskip\noindent
[LZ.1] T.P. Liu and K. Zumbrun,
{\it Nonlinear stability of an undercompressive shock for complex
Burgers equation,} Comm. Math. Phys. 168 (1995), no. 1, 163--186.
\medskip\noindent
[LZ.2] T.P. Liu and K. Zumbrun,
{\it On nonlinear stability of general undercompressive viscous shock waves,}
Comm.  Math. Phys. 174 (1995), no. 2, 319--345.
\medskip\noindent
[LZe.1] T.-P. Liu and Y. Zeng, 
{\it Large time behavior of solutions for general 
quasilinear hyperbolic--parabolic systems of conservation laws},
AMS memoirs 599 (1997).
\medskip\noindent
[LZe.2] T.-P. Liu and Y. Zeng, 
{\it Private communication} (manuscript in preparation).
\medskip\noindent
[M.1], A. Majda,
{\it The stability of multi-dimensional shock fronts -- a
new problem for linear hyperbolic equations,}
Mem. Amer. Math. Soc. 275 (1983).
\medskip\noindent
[M.2], A. Majda,
{\it The existence of multi-dimensional shock fronts,}
Mem. Amer. Math. Soc. 281 (1983).
\medskip\noindent
[M.3] A. Majda,
{\it Compressible fluid flow and systems of conservation laws in several
space variables,} Springer-Verlag, New York (1984), viii+ 159 pp.
\medskip\noindent
[MP] A. Majda and R. Pego,
{\it Stable viscosity matrices for systems of conservation laws},
J. Diff. Eqs. 56 (1985) 229--262.
\medskip\noindent
[MN] A. Matsumura and K. Nishihara, 
{\it On the stability of travelling wave solutions of a one-dimensional 
model system for compressible viscous gas,}
Japan J. Appl. Math. 2 (1985), no. 1, 17--25.
\medskip\noindent
[MZ.1] C. Mascia and K. Zumbrun,
{\it Pointwise Green's function bounds and stability of relaxation
shocks,} preprint (2001).
\medskip\noindent
[MZ.2] C. Mascia and K. Zumbrun,
{\it Pointwise Green's function bounds for shock profiles
with degenerate viscosity,} in preparation.
\medskip\noindent
[SZ] D. Serre and K. Zumbrun,
{\it Boundary layer stability in real-vanishing viscosity limit,}
to appear, Comm. Math. Phys. (2001).
\medskip\noindent
[SK] Y Shizuta and S. Kawashima,
{\it Systems of equations of hyperbolic--parabolic
type with applications to the discrete Boltzmann equation,}
Hokkaido Math. J. 14 (1984) 435--457.
\medskip\noindent
[Sm] J. Smoller, 
{\it Shock waves and reaction--diffusion equations,} 
Second edition, Grundlehren der Mathematischen Wissenschaften
[Fundamental Principles of Mathematical Sciences], 258. Springer-Verlag, 
New York, 1994. xxiv+632 pp. ISBN: 0-387-94259-9. 
\medskip\noindent
[Te] R. Temam, 
{\it Some developments on Navier-Stokes equations in the 
second half of the 20th century,} Development of mathematics
1950--2000, 1049--1106, Birkhäuser, Basel, 2000. 
\medskip\noindent
[Ze.1] Y. Zeng,
{\it $L^1$ asymptotic behavior of compressible, isentropic, 
viscous $1$-d flow,}
Comm. Pure Appl. Math. 47 (1994) 1053--1092.
\medskip\noindent
[ZH] K. Zumbrun and P. Howard,
{\it Pointwise semigroup methods and stability of viscous shock waves,}
Indiana Mathematics Journal V47 (1998) no. 4, 741--871.
\medskip\noindent
[Z.1] K. Zumbrun,  {\it Stability of viscous shock waves},
Lecture Notes, Indiana University (1998).
\medskip\noindent
[Z.2] K. Zumbrun, {\it Refined Wave--tracking and Nonlinear 
Stability of Viscous Lax Shocks}, to appear, Math. Anal. Appl. (2001);
preprint (1999).
\medskip\noindent
[Z.3] K. Zumbrun, {\it Multidimensional stability of
planar viscous shock waves}, TMR Summer School Lectures:
Kochel am See, May, 1999, to appear,
Birkhauser's Series: Progress in Nonlinear Differential
Equations and their Applications (2001), 207 pp.
\medskip\noindent
[Z.4] K. Zumbrun, {\it Stability of general undercompressive shocks
of viscous conservation laws,}
in preparation.
\medskip\noindent
[ZS] K. Zumbrun and D. Serre,
{\it Viscous and inviscid stability of multidimensional 
planar shock fronts,} Indiana Univ. Math. J. 48 (1999), no. 3,
937--992.
\endRefs
\enddocument